\documentclass[reqno]{amsart}
\usepackage{amssymb, latexsym}
\usepackage{hyperref}

\usepackage{backref}

\allowdisplaybreaks

\let\Re=\undefined\DeclareMathOperator*{\Re}{Re}
\let\Im=\undefined\DeclareMathOperator*{\Im}{Im}

\theoremstyle{plain}
\newtheorem{theorem}{Theorem}
\newtheorem{proposition}[theorem]{Proposition}
\newtheorem{lemma}[theorem]{Lemma}

\theoremstyle{definition}

\newtheorem{remark}[theorem]{Remark}

\newcounter{smalllist}

\numberwithin{equation}{section} \numberwithin{theorem}{section}
\begin{document}

\title[Wellposedness and singularities of the water wave equations]{Wellposedness and singularities of the water wave equations}
\author{Sijue Wu
}
\address{University of Michigan, Ann Arbor, MI}

\thanks {Financial support in part by NSF grants DMS-1101434, DMS-1361791}

\begin{abstract}

A class of water wave problems concerns the dynamics of the free interface separating an inviscid, incompressible and irrotational fluid, under the influence of gravity, from a zero-density region. In this  note, we present some recent methods and ideas developed concerning the local and global wellposedness of these problems, the focus is on the structural aspect of the equations.

\end{abstract}

\maketitle

\section{Introduction}

A class of water wave problems concerns the
motion of the 
interface separating an inviscid, incompressible, irrotational fluid,
under the influence of gravity, 
from a region of zero density (i.e. air) in 
$n$-dimensional space. It is assumed that the fluid region is below the
air region. Assume that
the density  of the fluid is $1$, the gravitational field is
$-{\bold k}$, where ${\bold k}$ is the unit vector pointing in the  upward
vertical direction, and at  
 time $t\ge 0$, the free interface is $\Sigma(t)$, and the fluid
occupies  region
$\Omega(t)$. When surface tension is
zero, the motion of the fluid is  described by 
\begin{equation}\label{euler}
\begin{cases}   \ \bold v_t + (\bold v\cdot \nabla) \bold v = -\bold k-\nabla P
\qquad  \text{on } \Omega(t),\ t\ge 0,
\\
\ \text{div}\,\bold v=0 , \qquad \text{curl}\,\bold v=0, \qquad  \text{on }
\Omega(t),\ t\ge 
0,
\\  
\ P=0, \qquad\qquad\qquad\qquad\qquad\text{on }
\Sigma(t) \\ 
\ (1, \bold v) \text{ 
is tangent to   
the free surface } (t, \Sigma(t)),
\end{cases}
\end{equation}
where $ \bold v$ is the fluid velocity, $P$ is the fluid
pressure. There is an important condition for these problems:
\begin{equation}\label{taylor}
-\frac{\partial P}{\partial\bold n}\ge 0
\end{equation}
pointwise on the interface, where $\bold n$ is the outward unit normal to the fluid interface 
$\Sigma(t)$ \cite{ta};
it is well known that when surface tension is neglected and the Taylor sign condition \eqref{taylor} fails, the water wave motion can be subject to the Taylor instability \cite{ ta, bi, bhl}.

The study on water waves dates back centuries. Early mathematical works include Stokes \cite{st}, Levi-Civita \cite{le}, and G.I. Taylor \cite{ta}.   Nalimov
\cite{na},  Yosihara \cite{yo} and Craig \cite{cr} proved local in time existence and uniqueness of solutions for the 2d water wave equation \eqref{euler} for small initial data. In  \cite{wu1, wu2}, we showed that for dimensions $n\ge 2$, the strong Taylor sign condition
 \begin{equation}\label{taylor-s}
-\frac{\partial P}{\partial\bold n}\ge c_0>0
\end{equation}
always holds for the
 infinite depth water wave problem \eqref{euler},  as long as the interface is non-self-intersecting and smooth; and the initial value problem of equation \eqref{euler} is locally well-posed in Sobolev spaces $H^s$, $s\ge 4$ for arbitrary given data.  Since then,
local wellposedness for water waves with additional effects such as the surface tension, bottom and non-zero vorticity,  under the assumption \eqref{taylor-s}\footnote{When there is surface tension, or bottom, or vorticity,  \eqref{taylor-s} does not hold, it needs to be assumed.} were obtained (c.f. \cite{am, cl, cs, ig1, la, li, ot, sz, zz}).  Alazard, Burq \& Zuily \cite{abz, abz14} proved local wellposedness of \eqref{euler} in low regularity Sobolev spaces where the interfaces are only in $C^{3/2}$. 
Recently, we proved the almost global and global wellposedness of \eqref{euler} for small, smooth and localized initial data for dimensions $n\ge 2$ \cite{wu3, wu4}; Germain, Masmudi \& Shatah obtained global existence for 3d water waves for a different class of small, smooth and localized data \cite{gms}.  Our 2d almost global existence result   has now been extended to global by Ionescu \& Pusateri  \cite{ip} and independently Alazard \& Delort  \cite{ad}, see \cite{hit, it} for an alternative proof. Finally, we  mention our most recent work on two dimensional water waves with angled crests \cite{wu5, kw, wu6}, in which we showed that for water waves with angled crests,  only the degenerate Taylor sign condition \eqref{taylor} holds, with degeneracy at the singularities on the interface. We proved an a priori estimate \cite{kw} and local existence \cite{wu6} for the 2d water wave equation \eqref{euler} in this framework. 

The advances of  water waves theory rely crucially on the understanding of the structure of the water wave equations. Indeed, \eqref{euler} is a nonlinear equation defined on moving domains, it is difficult to study it directly. A classical approach is to reduce from \eqref{euler} to an equation on the interface, and study solutions of the interface equation. Then use the incompressibility and irrotationality of the velocity field to recover the velocity in the fluid domain by solving a boundary value problem for the Laplace equation. However the fluid interface equation is itself a fully nonlinear and nonlocal equation, its structure is not easy to understood.  It is by achieving better understandings of the structure of this equation that has enabled us to apply analytical tools to deduce informations on the nature of the fluid motion. 

In this note we  describe the approach in  \cite{wu1, wu2, wu3, wu4}. Our focus is on the structural aspect of the work.
 It is clear that the 2d case is structurally simpler than 3d. 
 Our strategy  has been to first understand the two dimensional case, taking advantage  of complex analysis tools, in particular the Riemann mapping theorem, then use Clifford analysis to extend the 2d results to 3d. We note that although Riemann mapping is not available in 3d, it is by using it that has enabled us to understand the 2d case well enough to 
 develop an approach that extends to all dimensions. 
   
 We consider solutions of the water wave equation \eqref{euler} in the setting where
 $${\bold v}(\xi, t)\to 0,\qquad\text{as } |\xi|\to\infty$$
 and the interface $\Sigma(t)$ tends to the horizontal plane at infinity.\footnote{The problem with velocity $\bold v(\xi,t)\to (c,0)$ as $|\xi|\to\infty$ can be reduced to one with $\bold v$ tends to $0$ at infinity by the following observation: if $(\bold v, P)$ with $\Sigma(t): \xi=\xi(\cdot,t)$ is a solution of \eqref{euler}, then 
 $$\bold V=\bold v(\zeta+(c,0)t,t)-c, \quad \bold P=P(\zeta+(c,0)t,t),\quad\text{ with }\quad \Sigma(t)-(c,0)t: \zeta=\xi(\cdot, t)-(c,0)t$$
 is also a solution of \eqref{euler}.}
  
In section~\ref{section2},  we discuss the local wellposedness of \eqref{euler}, the focus is on deriving the quasilinear structure of equation \eqref{euler}, c.f. \cite{wu1, wu2, kw}. In section~\ref{section3}, we
consider the global in time behavior of solutions of \eqref{euler} in the regime of small waves, the focus is on understanding the nature of the nonlinearity of  equation \eqref{euler}, c.f. \cite{wu3, wu4}. We give some preparatory materials in the Appendices. In Appendix A, we  give some basic analysis tools such as estimates for commutators and operators involved in the equations. These inequalities tell us how a certain term behaves in terms of estimates, whether it is of higher order, or lower order etc., and guide us in our derivations of the structure of the equation. In Appendix B we  give some commutator identities that is used in our derivations. 

We use the following notations and conventions:  $[A, B]:=AB-BA$ is the commutator of operators $A$ and $B$; compositions are always in terms of the spatial variables and we write  for $f=f(\cdot, t)$, $g=g(\cdot, t)$, $f(g(\cdot,t),t):=f\circ g(\cdot, t):=U_gf(\cdot,t)$.

\section{Local wellposedness of the water wave equations}\label{section2}

From basic PDE theory we expect that in general, the Cauchy problem for a hyperbolic type PDE is locally solvable in Sobolev spaces, and we can solve it, for a short time, by energy estimates and an iterative argument; while the Cauchy problem for an elliptic type PDE is ill-posed. Hence in order to understand whether the water wave equation \eqref{euler} is uniquely solvable for a  positive time period for arbitrary given Cauchy data, it is crucial to understand its quasi-linear structure. In this section, we  derive the quasi-linear structure of the water wave equation \eqref{euler}, our focus is on the 2d case, since it is through a thorough understanding of this case that has enabled us to extend our work to 3d. We only give a brief description of how to extend the 2d derivations to 3d. We show  that the strong Taylor inequality \eqref{taylor-s}  always holds for $C^{1,\gamma}$ interfaces, while for singular interfaces only the weak Taylor inequality \eqref{taylor} holds; and this implies that the quasilinear structure of the water wave equation \eqref{euler} is of hyperbolic type in the regime of $C^{1,\gamma}$ interfaces and respectively of degenerate hyperbolic type in the regime that includes singular free surfaces.   The derivation given here is based on that in \cite{wu1, wu2}; due to the scope of this lecture note, we will only discuss the structural aspect of the work \cite{wu1, wu2}, and leave out the proof for the local wellposedness of \eqref{euler}. The interested reader may consult \cite{wu1, wu2} for the proof.


\subsection{The equation of the fluid interface in two space dimensions}\label{interface}
In two space dimensions, we identify $(x,y)$ with the complex number $x+iy$; $\Re z$, $\Im z$ are the real and imaginary parts of $z$; $\bar z=\Re z-i\Im z$ is the complex conjugate. 

Let the interface $\Sigma(t): z=z(\alpha, t)$, $\alpha\in\mathbb R$ be given by Lagrangian parameter $\alpha$, so $z_t(\alpha, t)={\bold v}(z(\alpha,t);t)$ is the velocity  of the fluid particles on the interface, $z_{tt}(\alpha,t)={\bold v_t + (\bold v\cdot \nabla) \bold v}(z(\alpha,t); t)$ is the acceleration; 
notice that $P=0$ on $\Sigma(t)$ implies that $\nabla P$ is normal to $\Sigma(t)$, therefore $\nabla P=-i\frak a z_\alpha$, where $\frak a =-\frac1{|z_\alpha|}\frac{\partial P}{\partial {\bold n}}$; we have from the first and third equation of \eqref{euler} that
\begin{equation}\label{interface-l}
z_{tt}+i=i\frak a z_\alpha
\end{equation}
The second equation of \eqref{euler}: $\text{div } \bold v=\text{curl } \bold v=0$ implies that $\bar {\bold v}$ is holomorphic in the fluid domain $\Omega(t)$; hence $\bar z_t$ is the boundary value of a holomorphic function in $\Omega(t)$.

Let $\Omega\subset \mathbb C$ be a domain with boundary $\Sigma: z=z(\alpha)$, $\alpha\in  I$, oriented clockwise. Let $\mathfrak H$ be the Hilbert transform associated to $\Omega$:
\begin{equation}\label{hilbert-t}
\frak H f(\alpha)=\frac1{\pi i}\, \text{pv.}\int\frac{z_\beta(\beta)}{z(\alpha)-z(\beta)}f(\beta)\,d\beta
\end{equation}
We have the following characterization of the trace of a holomorphic function on $\Omega$.

\begin{proposition}\cite{jour}\label{prop:hilbe}
a.  Let $g \in L^p$ for some $1<p <\infty$. 
  Then $g$ is the boundary value of a holomorphic function $G$ on $\Omega$ with $G(z)\to 0$ at infinity if and only if 
  \begin{equation}
    \label{eq:1571}
    (I-\mathfrak H) g = 0.
  \end{equation}

b. Let $ f \in L^p$ for some $1<p<\infty$. Then $ \mathbb P_H  f:=\frac12(I+\mathfrak H)  f$ is the boundary value of a holomorphic function $\frak G$ on $\Omega$, with $\frak G(z)\to 0$ as $|z|\to \infty$.

c. $\mathfrak H1=0$.
\end{proposition}
Observe  Proposition~\ref{prop:hilbe} gives that $\frak H^2=I$ on $L^p$. 

From Proposition~\ref{prop:hilbe}  the second equation of \eqref{euler} is equivalent to $\bar z_t=\mathfrak H {\bar z_t}$. Therefore the motion of the fluid interface $\Sigma(t): z=z(\alpha,t)$ is given by \footnote{
Equation \eqref{euler} and equation \eqref{interface-e} are equivalent, see \cite{wu1, wu2}. }
\begin{equation}\label{interface-e}
\begin{cases}
z_{tt}+i=i\frak a z_\alpha\\
\bar z_t=\frak H \bar z_t
\end{cases}
\end{equation}
\eqref{interface-e} is a fully nonlinear equation. To understand whether the equation is well posed, a usual strategy is to quasi-linearize the equation by differentiating. Notice that it can be hard to analyze the Hilbert transform $\frak H$ in the second equation since it  depends nonlinear nonlocally on the interface $z=z(\alpha,t)$,  this motivates us to use the Riemann mapping (c.f.\cite{wu1}).\footnote{In \cite{na, yo}, a quasilinear equation for the interface was derived for small and smooth waves in terms of the Lagrangian coordinates.  In \cite{wu1}, by using the Riemann mapping,  a more concise quasilinear equation was derived for all non-self-intersecting waves.  The approach described here is inspired by that in \cite{wu1},  it was used when we extended the work in \cite{wu1} to 3d \cite{wu2}. It is similar to that in \cite{wu1}, with the difference that in \cite{wu1} the derivation is in terms of the real component and here it is in both components.}

\subsection{Riemann mappings and the quasi-linear structure of 2d water waves}\label{riemann-mapping-q}
Let $\Phi(\cdot, t): \Omega(t)\to P_-$ be the Riemann mapping taking $\Omega(t)$ to the lower half plane $P_-$, satisfying $\lim_{z\to\infty}\Phi_z(z,t)=1$. Let 
$$h(\alpha;t):=\Phi(z(\alpha,t),t),$$
so $h:\mathbb R\to\mathbb R$ is a homeomorphism. Let $h^{-1}$ be defined by 
$$h(h^{-1}(\alpha',t),t)=\alpha',\quad \alpha'\in \mathbb R;$$
and 
$$Z(\alpha',t):=z\circ h^{-1}(\alpha',t),\quad Z_t(\alpha',t):=z_t\circ h^{-1}(\alpha',t),\quad Z_{tt}(\alpha',t):=z_{tt}\circ h^{-1}(\alpha',t)$$
be the reparametrization of the position, velocity and acceleration of the interface in the Riemann mapping variable $\alpha'$. Let
$$Z_{,\alpha'}(\alpha', t):=\partial_{\alpha'}Z(\alpha', t),\qquad Z_{t,\alpha'}(\alpha', t):=\partial_{\alpha'}Z_t(\alpha',t), \quad\text{etc.}$$
We note that $\Phi^{-1}(\alpha', t)=Z(\alpha', t)$. Notice that ${\bar {\bold v}}\circ \Phi^{-1}: P_-\to \mathbb C$ is holomorphic in the lower half plane $P_-$ with  ${\bar {\bold v}}\circ \Phi^{-1}(\alpha', t)={\bar Z}_t(\alpha',t)$. Precomposing \eqref{interface-l}  with $h^{-1}$ and applying Proposition~\ref{prop:hilbe} to $
{\bar {\bold v}}\circ \Phi^{-1}$ in $P_-$, we have the free surface equation in the Riemann mapping variable:
\begin{equation}\label{interface-r}
\begin{cases}
Z_{tt}+i=i\mathcal AZ_{,\alpha'}\\
\bar{Z}_t=\mathbb H \bar{Z}_t
\end{cases}
\end{equation}
where $\mathcal A\circ h=\frak a h_\alpha$ and $\mathbb H$ is the Hilbert transform associated with the lower half plane $P_-$:
$$\mathbb H f(\alpha')=\frac1{\pi i}\text{pv.}\int\frac1{\alpha'-\beta'}\,f(\beta')\,d\beta'.$$
\footnote{The advantage of using the Riemann mapping is that in $P_-$, the boundary value of a holomorphic function is characterized by $g=\mathbb Hg$, see Proposition~\ref{prop:hilbe}. The kernel of the Hilbert transform $\mathbb H$  is purely imaginary and is independent of the interface, it is easier to use $\mathbb H$ to understand the relations among various quantities and hence the quasi-linear structure of the free surface equation. 
}
As has been shown in \cite{wu1}, the quasi-linearization of \eqref{interface-e} or \eqref{interface-r} can be accomplished by just taking one derivative to $t$ to equation \eqref{interface-l}. 

Taking one derivative with respect to $t$ to \eqref{interface-l}, we get
\begin{equation}\label{quasi-l}
{\bar z}_{ttt}+i\frak a {\bar z}_{t\alpha}=-i\frak a_t {\bar z}_{\alpha}=\frac{\frak a_t}{\frak a} ({\bar z}_{tt}-i),
\end{equation}
the free surface equation is now \eqref{quasi-l} with the constraint $\bar z_t=\frak H \bar z_t$. Precomposing
with $h^{-1}$ on both sides of \eqref{quasi-l} we have in the Riemann mapping variable the free surface equation\footnote{\eqref{quasi-l} and \eqref{quasi-r} are equivalent to \eqref{interface-l} and \eqref{interface-r} provided the initial data for \eqref{quasi-l} and \eqref{quasi-r} stafisfy \eqref{interface-l} or \eqref{interface-r} at $t=0$.}
\begin{equation}\label{quasi-r}
\begin{cases}
{\bar Z}_{ttt}+i\mathcal A {\bar Z}_{t,\alpha'}=\frac{\frak a_t}{\frak a}\circ h^{-1} ({\bar Z}_{tt}-i)\\
\bar Z_t=\mathbb H \bar Z_t
\end{cases}
\end{equation}
We note that from the chain rule that for any function $f$,
$$U_{h^{-1}}\partial_t U_{h}f(\alpha', t)=(\partial_t+\mathcal B\partial_{\alpha'})f(\alpha', t)$$
where $\mathcal B=h_t\circ h^{-1}$. Hence  ${\bar Z}_{tt}=(\partial_t+\mathcal B\partial_{\alpha'}){\bar Z}_{t}$ and ${\bar Z}_{ttt}=(\partial_t+\mathcal B\partial_{\alpha'})^2{\bar Z}_{t}$. 
We will find $\mathcal B$, $\mathcal A$, and $\frac{\frak a_t}{\frak a}\circ h^{-1}$ in terms of ${\bar Z}_{t}$ and ${\bar Z}_{tt}$
and   show that \eqref{quasi-r} is a quasilinear equation for $\mathcal U=\bar Z_t$ with the right hand side consisting of lower order terms.\footnote{From Proposition~\ref{prop:hilbe}, $(I-\frak H)f$ measures how un-holomorphic  a function $f$ is.
The reason \eqref{quasi-l} is quasi-linear with the right hand side of lower order is that the two terms on the left 
hand side are "almost holomorphic" in the sense that $(I-\mathfrak H)(\bar z_{ttt}+i\frak a\partial_\alpha \bar z_t)$ is a commutator -- since by $(I-\frak H)\bar z_t=0$,  $(I-\mathfrak H)(\bar z_{ttt}+i\frak a\bar z_{t\alpha})=[\partial_t^2+i\frak a\partial_\alpha,\frak H]\bar z_t$, while the conjugate of the right hand side points  in the normal direction. We know a commutator is of lower order (c.f. Appendices A and B); and a holomorphic function with real part zero must be a constant; similarly an anti-holomorphic function with zero tangential part  on the boundary must be a constant. If the left hand side of \eqref{quasi-l} were holomorphic, then it would have to be zero. Now the left hand side is almost holomorphic with $(I-\frak H)(\bar z_{ttt}+i\frak a\partial_\alpha \bar z_t)$ a lower order term, then the right hand side must be a lower order term. We get this insight after our work in \cite{wu1}, not before.
} \footnote{We sometimes abuse notation and say a function $f$ is holomorphic if $f$ is the boundary value of a holomorphic function in $\Omega(t)$. 
}
To this end, we need some basic estimates for commutators, we leave these and some other preparatory materials in the Appendices. The reader may want to consult the Appendices before continuing.
\subsubsection{The quantity $\mathcal A$ and the Taylor sign condition} The basic idea of deriving the formulas for $\mathcal A$, $\mathcal B$ and $\frac{\frak a_t}{\frak a}\circ h^{-1}$ is to use the holomorphicity or almost holomorphicity of our quantities, and the  fact that $\Re (I-\mathbb H)f=f$ for real valued functions $f$. We will often write $(I-\mathbb H) (fg)$ as a commutator: $(I-\mathbb H) (fg)=[f, \mathbb H]g$ when $g$ satisfies $(I-\mathbb H)g=0$, since commutators are favorable in terms of estimates, see Appendix A.

Let $D_\alpha:=\frac 1{z_\alpha}\partial_\alpha $, and $D_{\alpha'}:=\frac 1{Z_{,\alpha'}}\partial_{\alpha'}$. Notice that for any holomorphic function $G$ in $\Omega(t)$ with boundary value $g(\alpha,t)=G(z(\alpha,t),t)$,
$$D_\alpha g(\alpha,t)=G_z(z(\alpha,t),t)$$
and for any function $f$, 
$$(D_\alpha f)\circ h^{-1}=D_{\alpha'} (f\circ h^{-1}).$$
 We note that $\mathcal A$ is related to the important quantity $-\frac{\partial P}{\partial {\bold n}}$ by $\mathcal A\circ h=\frak a h_\alpha$ and $\frak a=-\frac1{|z_\alpha|}\frac{\partial P}{\partial {\bold n}}$, therefore
\begin{equation}\label{taylor-f}
-\frac{\partial P}{\partial {\bold n}}\big |_{z=z(\alpha,t)}= (\mathcal A |Z_{,\alpha}|)\circ h.
\end{equation}
In this subsection we derive a formula for the quantity $\mathcal A$. This formula was first derived in \cite{wu1} to show that the strong Taylor sign condition \eqref{taylor-s} always holds for smooth non-self-intersecting interfaces.  It has also played a key role in our recent work on 2d water waves with angled crests \cite{kw, wu6}.

Taking complex conjugate of the first equation in \eqref{interface-r} then multiply by $Z_{,\alpha'}$ yields
\begin{equation}\label{interface-a1}
Z_{,\alpha'}({\bar Z}_{tt}-i)=-i\mathcal A|Z_{,\alpha'}|^2:=-i A_1.
\end{equation}
The left hand side of \eqref{interface-a1} is almost holomorphic since $Z_{,\alpha'}$ is the boundary value of the holomorphic function $(\Phi^{-1})_{z'}$ and $\bar z_{tt}$ is the time derivative of the holomorphic function $\bar z_t$. We explore the almost holomorphicity of $\bar z_{tt}$ by expanding.  Let $F=\bar {\bold v}$, we know $F$ is holomorphic in $\Omega(t)$, and $\bar z_t=F(z(\alpha, t),t)$, so
\begin{equation}\label{eq:1}
\bar z_{tt}=F_t(z(\alpha, t),t)+F_z(z(\alpha, t),t) z_t(\alpha, t),\qquad \bar z_{t\alpha}=F_z(z(\alpha, t),t) z_\alpha(\alpha, t)
\end{equation}
therefore \begin{equation}\label{eq:2}
\bar z_{tt}= F_t\circ z+ \frac{\bar z_{t\alpha}}{z_\alpha} z_t.
\end{equation}
 Precomposing with $h^{-1}$, subtracting $-i$, then multiplying by $Z_{,\alpha'}$, we have
$$Z_{,\alpha'}({\bar Z}_{tt}-i)= Z_{,\alpha'} F_t\circ Z+ Z_t {\bar Z}_{t,\alpha'}-i Z_{,\alpha'}=-iA_1 $$
Apply $(I-\mathbb H)$ to both sides of the equation. Notice that  $F_t\circ Z$ is the boundary value of the holomorphic function  $F_t\circ \Phi^{-1}$,  so $(I-\mathbb H)(Z_{,\alpha'} F_t\circ Z)=0$, $ (I-\mathbb H)Z_{,\alpha'}=1$,\footnote{We know $(\Phi^{-1})_{z'}\to 1$ as $z'\to \infty$; we assume a priori that $(\Phi^{-1})_{z'}F_t\circ \Phi^{-1}\to 0$ as $z'\to \infty$. As proved in \cite{wu1} equation \eqref{interface-r} is well posed in this regime.}  therefore
$$-i(I-\mathbb H) A_1= (I-\mathbb H)(Z_t {\bar Z}_{t,\alpha'})-i$$
Taking imaginary parts on both sides and using the fact $(I-\mathbb H){\bar Z}_{t,\alpha'}=0$ \footnote{Because $(I-\mathbb H){\bar Z}_{t}=0$.}  to rewrite $(I-\mathbb H)(Z_t {\bar Z}_{t,\alpha'})$ as $[Z_t,\mathbb H]{\bar Z}_{t,\alpha'}$ yields
\begin{equation}\label{A_1}
A_1=1-\Im [Z_t,\mathbb H]{\bar Z}_{t,\alpha'}
\end{equation}
Compute
\begin{equation}
\begin{aligned}
-\Im[Z_t,\mathbb H]{\bar Z}_{t,\alpha'}&=-\frac1{\pi} \Re\int\frac{(Z_t(\alpha',t)-Z_t(\beta',t))\partial_{\beta'}({\bar Z}_t(\alpha',t)-{\bar Z}_t(\beta',t))}{\alpha'-\beta'}\,d\beta\\&
=-\frac1{2\pi} \int\frac{\partial_{\beta'}|{ Z}_t(\alpha',t)-{ Z}_t(\beta',t)|^2}{\alpha'-\beta'}\,d\beta\\
&=\frac1{2\pi }\int \frac{|Z_t(\alpha', t)-Z_t(\beta', t)|^2}{(\alpha'-\beta')^2}\,d\beta'
\end{aligned}
\end{equation}
where in the last step we used integration by parts. We conclude 
\begin{proposition}(c.f. \cite{wu1}, Lemma 3.1)\label{prop:a1} We have
1. \begin{equation}\label{a1}
A_1=1-\Im [Z_t,\mathbb H]{\bar Z}_{t,\alpha'}=1+\frac1{2\pi }\int \frac{|Z_t(\alpha', t)-Z_t(\beta', t)|^2}{(\alpha'-\beta')^2}\,d\beta'\ge 1.
\end{equation}
2. \begin{equation}\label{taylor-formula}
-\frac{\partial P}{\partial\bold n}= \frac{A_1}{|Z_{,\alpha'}|};
\end{equation}
in particular if the interface $\Sigma(t)\in C^{1,\gamma}$ for some $\gamma>0$, then the strong Taylor sign condition \eqref{taylor-s} holds.
\end{proposition}
From \eqref{interface-a1} we have
\begin{equation}\label{a}
\mathcal A=\frac{A_1}{|Z_{,\alpha'}|^2}=\frac{|{\bar Z}_{tt}-i|^2}{A_1}
\end{equation}
with $A_1$ given by \eqref{a1}. 
\subsubsection{The quantity $\mathcal B=h_t\circ h^{-1}$} The quantity $\mathcal B=h_t\circ h^{-1}$ can be calculated similarly.
Recall $h(\alpha,t)=\Phi(z(\alpha,t),t)$, so 
$$h_t=\Phi_t\circ z+(\Phi_z\circ z) z_t,\qquad h_\alpha=(\Phi_z\circ z) z_\alpha$$
hence
$h_t= \Phi_t\circ z+ \frac{h_\alpha}{z_\alpha}  z_t$. Precomposing with $h^{-1}$ yields
\begin{equation}\label{b1}
h_t\circ h^{-1}=\Phi_t\circ Z+ \frac{Z_t}{Z_{,\alpha'}}.
\end{equation}
Now $(I-\mathbb H)\Phi_t\circ Z=(I-\mathbb H)(\frac 1{Z_{,\alpha'}}-1)=0$ since $\Phi_t\circ Z$ and $\frac 1{Z_{,\alpha'}}-1$ are boundary values of the holomorphic functions $\Phi_t\circ \Phi^{-1}$ and $\frac1{(\Phi^{-1})_{z'}}-1$ respectively.\footnote{Again, we assume a priori that $\Phi_t\circ \Phi^{-1}(z',t)\to 0$ as $z'\to\infty$. It has been proved in \cite{wu1} that such solutions exist.} Apply $(I-\mathbb H)$ to both sides of \eqref{b1} then take the real parts; rewriting $ (I-\mathbb H)(Z_t(\frac{1}{Z_{,\alpha'}}-1))$ as $[Z_t,\mathbb H](\frac1{Z_{,\alpha'}}-1)$,  we get
\begin{equation}\label{b}
\mathcal B=h_t\circ h^{-1}=\Re \big([Z_t,\mathbb H](\frac1{Z_{,\alpha'}}-1)\big)+2\Re Z_t.
\end{equation}
Here  we used the fact that $\Re(I-\mathbb H)Z_t=2\Re Z_t$, since $(I+\mathbb H)Z_t=\overline{(I-\mathbb H)\bar Z_t}=0$.
\subsubsection{The quantity $\frac{\frak a_t}{\frak a}\circ h^{-1}$} We analyze $\frac{\frak a_t}{\frak a}\circ h^{-1}$ similarly. Start with \eqref{quasi-l} or the first equation of \eqref{quasi-r}, and expand the left hand side using $\bar z_t=F(z(\alpha,t), t)$. Taking one more derivative to $t$ to the first equation in \eqref{eq:1}, we get
$${\bar z}_{ttt}=(F_{zz}\circ z) z_t^2+2(F_{tz}\circ z)z_t+(F_z\circ z)z_{tt}+F_{tt}\circ z.$$
From \eqref{eq:2}, the second equation of \eqref{eq:1}, and the holomorphicity of $F_z$, $F_t$, we have
$$F_z\circ z=D_\alpha \bar z_t, \quad F_{zz}\circ z= D_\alpha^2 \bar z_t,\quad\text{and } F_{tz}\circ z=D_\alpha({\bar z}_{tt}-(D_\alpha \bar z_t)z_t)$$
therefore
$${\bar z}_{ttt}=(D_\alpha^2 \bar z_t) z_t^2+2z_tD_\alpha({\bar z}_{tt}-(D_\alpha \bar z_t)z_t)+(D_\alpha \bar z_t)z_{tt}+F_{tt}\circ z.$$
Precomposing with $h^{-1}$ yields
\begin{equation}\label{eq:3}
{\bar Z}_{ttt}=(D_{\alpha'}^2 \bar Z_t) Z_t^2+2Z_tD_{\alpha'}({\bar Z}_{tt}-(D_{\alpha'} \bar Z_t)Z_t)+(D_{\alpha'} \bar Z_t)Z_{tt}+F_{tt}\circ Z.
\end{equation}
Now multiply the first equation of \eqref{quasi-r} by $Z_{,\alpha'}$, then substitute in \eqref{eq:3}. We have, by \eqref{interface-a1},
\begin{equation}\label{eq:4}
\begin{aligned}
Z_{,\alpha'}\{(D_{\alpha'}^2 \bar Z_t) Z_t^2+2Z_tD_{\alpha'}({\bar Z}_{tt}-(D_{\alpha'} \bar Z_t)Z_t)+(D_{\alpha'} \bar Z_t)Z_{tt}+&F_{tt}\circ Z+i\mathcal A\bar Z_{t,\alpha'}\}\\&=-i A_1 \frac{\frak a_t}{\frak a}\circ h^{-1}.
\end{aligned}
\end{equation}
Apply $(I-\mathbb H)$ to both sides of \eqref{eq:4}. Using the fact that $(I-\mathbb H)(\partial_{\alpha'}D_{\alpha'} \bar Z_t)=0$, $(I-\mathbb H)\partial_{\alpha'}({\bar Z}_{tt}-(D_{\alpha'} \bar Z_t)Z_t)=0$, $(I-\mathbb H) \partial_{\alpha'} \bar Z_t=0$ and $(I-\mathbb H)F_{tt}\circ Z=0$ because of the holomorphicity of $\partial_{z'}(F_z\circ \Phi^{-1})$, $\partial_{z'}(F_t\circ \Phi^{-1})$, $\partial_{z'}(F\circ \Phi^{-1})$ and $F_{tt}\circ\Phi^{-1}$, and the identity $Z_{tt}+i=i\mathcal AZ_{,\alpha'}$ \eqref{interface-r}, we rewrite each term on the left as commutator and get
\begin{equation}\label{eq:5}
\begin{aligned}
 \ [Z_t^2,\mathbb H]\partial_{\alpha'}D_{\alpha'} \bar Z_t+2[Z_t,\mathbb H]\partial_{\alpha'}({\bar Z}_{tt}-(D_{\alpha'} \bar Z_t)Z_t)&+2[Z_{tt},\mathbb H]\partial_{\alpha'} \bar Z_t\\&=(I-\mathbb H)(-i A_1 \frac{\frak a_t}{\frak a}\circ h^{-1}).
\end{aligned}
\end{equation}
For the sake of estimates we need to further rewrite the first term and the second part of the second term. Using integration by parts on each term, after cancelations we have
\begin{align*}
[Z_t^2,\mathbb H]\partial_{\alpha'}D_{\alpha'} \bar Z_t-2[Z_t,\mathbb H]\partial_{\alpha'}((D_{\alpha'} \bar Z_t)Z_t)&=-\frac1{\pi i}\int\frac{(Z_t(\alpha',t)-Z_t(\beta',t))^2}{(\alpha'-\beta')^2} D_{\beta'} \bar Z_t(\beta',t)\,d\beta'\\&:=-[Z_t, Z_t; D_{\alpha'} \bar Z_t].
\end{align*}
Therefore we can write \eqref{eq:5} as
\begin{equation}\label{eq:6}
2[Z_t,\mathbb H]{\bar Z}_{tt,\alpha'}+2[Z_{tt},\mathbb H]\partial_{\alpha'} \bar Z_t-
[Z_t, Z_t; D_{\alpha'} \bar Z_t]
=(I-\mathbb H)(-i A_1 \frac{\frak a_t}{\frak a}\circ h^{-1}).
\end{equation}
Taking imaginary parts on both sides and dividing by $-A_1$ yields
\begin{equation}\label{at}
\frac{\frak a_t}{\frak a}\circ h^{-1}=\frac{-\Im( 2[Z_t,\mathbb H]{\bar Z}_{tt,\alpha'}+2[Z_{tt},\mathbb H]\partial_{\alpha'} \bar Z_t-
[Z_t, Z_t; D_{\alpha'} \bar Z_t])}{A_1}.
\end{equation}
\subsubsection{The quasilinear equation} Sum up \eqref{quasi-r}, \eqref{a1}, \eqref{a}, \eqref{b} and \eqref{at}, we have the equation for the free interface: 
\begin{equation}\label{quasi-linear}
\begin{cases}
(\partial_t+\mathcal B\partial_{\alpha'})^2\bar Z_t+ i\mathcal A\partial_{\alpha'} \bar Z_t=g\\
\bar Z_t=\mathbb H\bar Z_t
\end{cases}
\end{equation}
where
\begin{equation}\label{aux-r}
\begin{cases}
\mathcal A=\frac{|{\bar Z}_{tt}-i|^2}{A_1},\quad Z_{tt}=(\partial_t+\mathcal B\partial_{\alpha'}) Z_t\\
A_1=1-\Im [Z_t,\mathbb H]{\bar Z}_{t,\alpha'}=1+\frac1{2\pi }\int \frac{|Z_t(\alpha', t)-Z_t(\beta', t)|^2}{(\alpha'-\beta')^2}\,d\beta'\\
\mathcal B=h_t\circ h^{-1}=\Re\big([Z_t,\mathbb H](\frac1{Z_{,\alpha'}}-1)\big)+2\Re Z_t\\
\frac1{Z_{,\alpha'}}=i\frac{\bar Z_{tt}-i}{A_1}\\
g=(\bar Z_{tt}-i) \frac{-\Im( 2[Z_t,\mathbb H]{\bar Z}_{tt,\alpha'}+2[Z_{tt},\mathbb H]\partial_{\alpha'} \bar Z_t-
[Z_t, Z_t; D_{\alpha'} \bar Z_t])}{A_1}.
\end{cases}
\end{equation}
 From the inequalities in Appendix A, we know in the regime where the interface is $C^{1,\gamma}$, $\gamma>0$, \eqref{quasi-linear}-\eqref{aux-r} is a quasilinear equation of the conjugate velocity $\bar Z_t$,  with the right hand side of \eqref{quasi-linear} consisting of  lower order terms.  \eqref{quasi-linear} is of hyperbolic type since $$\mathcal A=\frac{A_1}{|Z_{,\alpha'}|^2}\ge c_0>0$$
for $C^{1,\gamma}$ interfaces, and $i\partial_{\alpha'}\bar Z_t=|\partial_{\alpha'}|\bar Z_t$.\footnote{This is because $i\partial_{\alpha'}\bar Z_t=i\partial_{\alpha'}\mathbb H \bar Z_t=|\partial_{\alpha'}|\bar Z_t$.}
Local well-posedness of \eqref{quasi-linear}-\eqref{aux-r} for $(\bar Z_t,\bar Z_{tt})\in C([0, T], H^{s+1/2}\times H^s)$, $s\ge 4$ has been proved by the energy method and an iterative argument, we refer the reader to \cite{wu1} for details. 

We make the following remarks concerning some recent works.
\begin{remark}
1. The difference between the quasilinear equation \eqref{quasi-linear}-\eqref{aux-r} and the quasilinear equation (4.6)-(4.7) of \cite{wu1} is that \eqref{quasi-linear}-\eqref{aux-r} is in terms of both components of $Z_t$ and (4.6)-(4.7) of \cite{wu1} is in terms of the real component of $Z_t$. (4.6)-(4.7) of \cite{wu1} has been written in
a way so that it is easy to prove its equivalence with the interface equation \eqref{interface-e} or equation (1.7)-(1.8) of \cite{wu1}; see \S6 of \cite{wu1}.

2. \eqref{quasi-linear}-\eqref{aux-r} is an equation for the conjugate velocity $\bar Z_t$ and conjugate acceleration $\bar Z_{tt}$, the interface doesn't appear explicitly, so a solution of  \eqref{quasi-linear}-\eqref{aux-r} can exist even when $Z=Z(\cdot, t)$ becomes self-intersecting.\footnote{$Z=Z(\cdot, t)$ is defined by $z(\cdot, t)=z(\cdot,0)+\int_0^t z_s(\cdot, s)\,ds$, where $z_t=Z_t\circ h$, $Z=z\circ h^{-1}$; and $h_t=\mathcal B(h, t)$.}  Checking through the derivation above we see that we arrived at 
\eqref{quasi-linear}-\eqref{aux-r} from \eqref{euler} using only the following properties of the domain: 1. there is a conformal mapping taking the fluid region $\Omega(t)$ to $P_-$; 2. $P=0$ on $\Sigma(t)$. We note that  $z\to z^{1/2}$ is a conformal map that takes the region $\mathbb C\setminus \{z=x+i0, x>0\}$ to the upper half plane; a domain with its boundary self-intersecting at the positive real axis can therefore be mapped conformally onto the lower half plane $P_-$. Taking such a domain as the initial fluid domain, 
assuming  $P=0$ on $\Sigma(t)$ even when $\Sigma(t)$ self-intersects,\footnote{We note that when $\Sigma(t)$ self-intersects, the condition $P=0$ on $\Sigma(t)$ is unphysical.} one can still solve equation \eqref{quasi-linear}-\eqref{aux-r} for a short time. Indeed this is the main idea in the work of \cite{cf}. Using this idea and the time reversibility of the water wave equation, by choosing an appropriate initial velocity field that pulls the initial domain apart,  \cite{cf} proved the existence of "splash" and "splat" singularities starting from a smooth non-self-intersecting fluid region. 

3. The above derivation applies  to fluid domains with arbitrary non-self-intersecting boundaries. We have from \eqref{taylor-formula} and \eqref{a1} that the Taylor sign condition \eqref{taylor} always holds, as long as the interface is non-self-intersecting. Assume that $\Sigma(t)$ is non-self-intersecting with angled crests, assume the interior angle at a crest is $\nu$. We know the Riemann mapping $\Phi^{-1}$ (we move the singular point to the origin) behaves like 
$$\Phi^{-1}(z')\approx (z')^r,\qquad \text{with } \nu=r\pi$$
near the crest, so $Z_{,\alpha'}\approx (\alpha')^{r-1}$ near the crest. From \eqref{interface-a1} and the fact $A_1\ge 1$, we 
can conclude that the interior angle at the crest must be $\le \pi$ if the acceleration $|Z_{tt}|\ne\infty$;
we can also conclude that $-\frac{\partial P}{\partial\bold n}=0$ at the singularities where the interior angles are $<\pi$, therefore in the regime that includes singular free surfaces, the quasilinear equation \eqref{quasi-linear}-\eqref{aux-r} is degenerate hyperbolic, c.f. \cite{kw, wu6}. In \cite{kw, wu6} we proved an a priori estimate and  the local existence for 2d water waves in the regime including interfaces with angled crests, showing that the water wave equation \eqref{euler} admit such solutions.

4. The quasi-linear equation for water waves in the periodic setting can be derived similarly, see \cite{kw}.

 5. The Riemann mapping variable is used in recent work \cite{hit, it}.

\end{remark}

\subsection{The quasilinear equation in Lagrangian coordinates}\label{quasilinear-l}
 In order to extend our work for 2d water waves to 3d, we  need a derivation that does not rely on the Riemann mapping. Upon checking the derivation in \S\ref{riemann-mapping-q}, we see that all we have done is to apply $\Re(I-\mathbb H)$ to calculate the parameters in the equations. We certainly can do the same calculations with $(I-\frak H)$, Riemann mapping is not needed. 
  
  We now use this idea to analyze $\frak a_t|z_\alpha|$ and show that indeed \eqref{quasi-l} is quasilinear with the right hand side of lower order. 
 
 First by \eqref{interface-l} and Proposition~\ref{prop:a1} we have 
 \begin{equation}\label{al}
 \frak a{|z_\alpha|}={|z_{tt}+i|},\qquad \text{and}\quad -i\frac{\bar z_\alpha}{|z_\alpha|}=\frac{\bar z_{tt}-i}{|z_{tt}+i|}.
 \end{equation}
 We apply $(I-\frak H)$ to \eqref{quasi-l}. Using 
$\bar z_t= \frak {H} \bar z_t$ and Proposition~\ref{lemma1} in Appendix B,  we have
\begin{equation}\label{eq:11}
\begin{aligned}
 & (I-\frak{H})(-i{\frak a}_t \bar{z}_\alpha)=(I-\frak{H})(\bar{z}_{ttt}+ i \frak{a} \bar{z}_{t\alpha})
\\&=[\partial_t^2+ i\frak{a} \partial_\alpha, \frak{H}]\bar z_t\\&=
2[z_{tt}, \frak {H}]\frac{\bar z_{t\alpha}}{z_\alpha}+2[z_t, \frak {H}]\frac{\bar z_{tt\alpha}}{z_\alpha} -\frac 1{\pi  i}\int(\frac{z_t(\alpha, t)-z_t(\beta, t)}{z(\alpha,t)-z(\beta,t)})^2\bar z_{t\beta}\, d\beta\end{aligned} 
\end{equation}
Multiply both sides of \eqref{eq:11} by $i\frac{z_\alpha}{|z_\alpha|}$ and take the real parts. 
Since $\frak a$ and $\frak a_t$ are real valued, we have
\begin{equation}\label{2dat}
\begin{aligned}
(I+&\frak K^*)( \frak a_t | z_\alpha|)=\\
& {\Re}(\frac{i z_\alpha}{|z_\alpha|}\{2[z_{tt}, \frak {H}]\frac{\bar z_{t\alpha}}{z_\alpha}+2[z_t, \frak {H}]\frac{\bar z_{tt\alpha}}{z_\alpha} -\frac 1{\pi  i}\int(\frac{z_t(\alpha, t)-z_t(\beta, t)}{z(\alpha,t)-z(\beta,t)})^2\bar z_{t\beta}\, d\beta\})
\end{aligned}
\end{equation}
where
$$\frak K^*f(\alpha,t)=p.v.\int {\Re}\{\frac{-1}{\pi i} \frac{z_\alpha}{|z_\alpha|}\frac{|z_\beta(\beta,t)|}{(z(\alpha,t)-z(\beta,t))}\}f(\beta,t)\,d\beta$$ 
is the adjoint of the double layer potential operator $\frak K$ in $L^2(\Sigma(t), dS)$. We know $I+\frak K^*$ is invertible on $L^2(\Sigma(t), dS)$ (cf. \cite{fol, ke}).   The second equation of \eqref{al}, \eqref{2dat} and the estimates in Appendix A show that  $\frak{a}_t|z_\alpha|$ has the same regularity as that of  $\bar z_{tt}$ and $\bar z_t$.

 We rewrite \eqref{quasi-l} with the constraint $\bar z_t=\frak H\bar z_t$ as 
\begin{equation}\label{2dquasilinear}
\begin{cases} \bar{z}_{ttt}+ i \frak{a} \bar {z}_{t\alpha}
= \frac{ \bar z_{tt}-i}{| z_{tt}+i|}\frak a_t | z_\alpha|\\
\bar z_t= \frak {H} \bar z_t,\end{cases}
\end{equation}
where $\frak a |z_\alpha|$ is given by \eqref{al}  and  $\frak a_t | z_\alpha|$ is given by \eqref{2dat}. Since $\bar z_t$ is holomorphic, $i\frac{1}{|z_\alpha|}\partial_\alpha \bar z_t=\nabla_{\bold n}\bar z_t$. By Green's identity, the Dirichlet-Neumann operator $\nabla_{\bold n}$ is a positive operator. From Proposition~\ref{prop:a1}  we know $\frak {a}{|z_\alpha|}=-\frac{\partial P}{\partial{\bold n}}\ge c_0>0$ in the regime of $C^{1,\gamma}$ interfaces.
Therefore \eqref{2dquasilinear} is a quasilinear system of hyperbolic type, with the right hand side of the first equation in \eqref{2dquasilinear} consisting of terms of lower order derivatives of $ \bar z_t$.  
The local in time wellposedness of   \eqref{2dquasilinear}-\eqref{2dat}-\eqref{al}  in Sobolev spaces (with $(z_t, z_{tt})\in C([0,T], H^{s+1/2}\times H^s)$, $s\ge 4$ ) can then be proved
 by energy estimates and a fixed point iteration argument. 
\footnote{A proof of the local wellposedness of the 3d counterpart of \eqref{2dquasilinear}-\eqref{2dat}-\eqref{al} is carried out in \cite{wu2}.} 

This derivation has been extended to 3d using Clifford analysis, cf. \cite{wu2}, we will give a brief discussion on how to do this in \S\ref{3dquasi-l}. 

Before ending this subsection, we mention that the quasilinear system \eqref{2dquasilinear}-\eqref{2dat}-\eqref{al} is coordinate invariant. 

For fixed $t$, let $k=k(\alpha, t): \mathbb R\to \mathbb R$ be a diffeomorphism with $k_\alpha >0$.  Let $k^{-1}$ be such that $k\circ k^{-1}(\alpha, t)=\alpha$. Define 
\begin{equation}\label{2dzeta}
\zeta:=z\circ k^{-1}, \qquad b:=k_t\circ k^{-1} \quad\text{and}\quad A\circ k:={\frak a}k_\alpha.
\end{equation}  
Let \begin{equation}\label{2dzetaa1}
D_t:=U_k^{-1}\partial_tU_k:=\partial_t+b\partial_\alpha
\end{equation}
 be the material derivative. By a simple  application of the chain rule, we have
$$U_k^{-1}(\partial_t^2+i\frak{a}\partial_\alpha)U_k=D_t^2+iA\partial_\alpha,$$
and equation \eqref{2dquasilinear}  becomes
\begin{equation}\label{k2dquasilinear}
\begin{cases} (D_t^2+i A \partial_\alpha) \overline{D_t\zeta}=  ({\frak a_t | z_\alpha|})\circ k^{-1}\dfrac{\overline{D_t^2\zeta}-i}{|D_t^2\zeta+i|}\\
\overline {D_t\zeta}={\mathcal H}\overline {D_t\zeta}\end{cases}
\end{equation}
with
\begin{equation}\label{k2dat}
\begin{aligned} 
(I+&\mathcal K^*)( (\frak a_t | z_\alpha|)\circ k^{-1})= \\&
{Re}(\frac{i \zeta_\alpha}{|\zeta_\alpha|}\{2[ D_t^2\zeta, {\mathcal H}]\frac{ \partial_{\alpha} \overline {D_t\zeta}}{\zeta_\alpha}+2[ D_t\zeta, {\mathcal H}]\frac{\partial_{\alpha}\overline {D_t^2\zeta}}{\zeta_\alpha}\\& -\frac 1{\pi  i}\int\big(\frac{ D_t\zeta(\alpha, t)- D_t\zeta(\beta, t)}{\zeta(\alpha,t)-\zeta(\beta,t)}\big)^2 \partial_{\beta}\overline {D_t\zeta}(\beta,t)\, d\beta\})
\end{aligned}
\end{equation}
and 
\begin{equation}\label{khilbert}
{\mathcal H} f(\alpha,t)=U_k^{-1}\frak H U_kf(\alpha,t)=\frac 1{\pi i}p.v.\int\frac {f(\beta, t) \zeta_\beta(\beta,t)}{\zeta(\alpha,t)-\zeta(\beta, t)}\,d\,\beta,
\end{equation}
 \begin{equation}
 \mathcal K^*f(\alpha,t)=p.v.\int {Re}\{\frac{-1}{\pi i} \frac{\zeta_\alpha}{|\zeta_\alpha|}\frac{|\zeta_\beta(\beta,t)|}{(\zeta(\alpha,t)-\zeta(\beta,t))}\}f(\beta,t)\,d\beta.
 \end{equation}
  Notice the remarkable similarities between equations \eqref{2dquasilinear}-\eqref{2dat} and \eqref{k2dquasilinear}-\eqref{k2dat}. In particular,  the structures of the terms in \eqref{2dquasilinear}-\eqref{2dat} do not change under the change of variables. This makes it  convenient for us to  work in another coordinate system and to choose a different coordinate system when there is advantage to do so. In fact, this has been used in our study of the global in time behavior of water waves \cite{wu3, wu4}. 


\subsection{The quasi-linear equation for  3d water waves}\label{3dquasi-l}
We derive the quasi-linear equation for 3d water waves by carrying out the same procedure as for 2d. 
We  first need to write down the 3d counterpart of the interface equation \eqref{interface-e}. While equation \eqref{interface-l} is readily available in 3d, to write down the second equation, we need a suitable counterpart in 3d of the equation for the trace on the interface of the velocity field $\bold v$ that satisfies $\text{div\,}\bold v=0$ and $\text{curl\,} \bold v=0$.  This leads us to  Clifford analysis.



Let's recall the basics of Clifford algebra, or in other words, the algebra of quaternions  
$\mathcal C(V_2)$ (c.f. \cite{gm}). Let $\{1, e_1, e_2, e_3\}$ be the basis of $\mathcal C(V_2)$, satisfying 
\begin{equation}\label{product}
e_i^2=-1,\quad\text{for }i=1,2,3, \qquad e_i e_j =-e_j e_i,\quad i\ne j, \qquad e_3=e_1e_2.
\end{equation}
Let $\mathcal D= \partial_x e_1+ \partial_y e_2+ \partial_z e_3$. By definition, a Clifford-valued function 
$F: \Omega\subset \mathbb R^3\to \mathcal C(V_2)$ is Clifford analytic in domain $\Omega$ iff $\mathcal DF=0$ in $\Omega$. Therefore 
$F=\sum_{i=1}^3 f_i e_i$ is Clifford analytic in $\Omega$ if and only if $\text{div} F=0$ and $\text{curl} F=0$ in $\Omega$. Furthermore a function $F$ is the trace of a 
Clifford analytic function in $\Omega$ if and only if $F=\frak H_{\Sigma} F$, where
\begin{equation}\label{3dhilbert}
\frak H_{\Sigma} g(\alpha,\beta )=
p.v. \iint K(\eta(\alpha',\beta') -\eta(\alpha,  \beta))\,(\eta'_{\alpha'}\times\eta'_{\beta'})\,g( \alpha',\beta'
)\,d\alpha'd\beta'
\end{equation}
is the 3d version of the Hilbert transform on $\Sigma=\partial\Omega: \eta=\eta(\alpha,\beta), \ (\alpha,\beta)\in \mathbb R^2$,  with  normal $\eta_{\alpha}\times\eta_{\beta}$ pointing out of $\Omega$, and 
$$\Gamma(\eta)=-\frac1{\omega_3 |\eta|}, \qquad K(\eta)=-2\mathcal D \Gamma(\eta)=-\frac{2}{\omega_3}\frac{\eta}{|\eta|^3},$$ 
 $\omega_3$ is the surface area of the unit sphere in $\mathbb R^3$.\footnote{\eqref{3dhilbert} is similar to the 2d Hilbert transform \eqref{hilbert-t} in the sense that in 2d, the fundamental solution for Laplace equation is $\Gamma(z)=\frac1{2\pi}\ln{|z|}$, $K(z)=-2(\partial_x-i\partial_y)\big(\frac1{2\pi}\ln{|z|}\big)=-\frac1{\pi z}$, the outward normal is $iz_\alpha$, and $K(z(\alpha)-z(\beta))(iz_\beta)=\frac1{\pi i}\frac{iz_\beta}{z(\alpha)-z(\beta)}$ is the kernel of the Hilbert transform \eqref{hilbert-t}.
}
As in the 2d case, if $\xi=\xi(\alpha,\beta, t)$, $(\alpha,\beta)\in \mathbb R^2$ is the free interface $\Sigma(t)$ in Lagrangian coordinates $(\alpha,\beta)$ at time $t$, with $N=\xi_\alpha\times\xi_\beta$ pointing out of the fluid domain, we can rewrite the 3D water wave system \eqref{euler} ($n=3$) as 
\begin{equation}\label{3dww1}
\begin{cases}
\xi_{tt}+e_3=\frak a N\\
\xi_t=\frak H_{\Sigma(t)}\xi_t
\end{cases}
\end{equation}
where $\frak a=-\frac1{|N|}\frac{\partial P}{\partial \text{\bf n}}$.

Differentiating the first equation with respect to $t$ yields
\begin{equation}\label{3dquasi-linear}
\begin{cases}
\xi_{ttt}-\frak a N_t=\frak a_t N\\
\xi_t=\frak H_{\Sigma(t)}\xi_t.
\end{cases}
\end{equation}
This is the 3d counterpart of the 2d quasilinear equation \eqref{quasi-l} with constraint $\bar z_t=\frak H\bar z_t$.   It has been proved in \cite{wu2} that $ N_t=-|N|\nabla_{\bold n}\xi_t$,\footnote{In 2d, we know $iz_\alpha$ is a normal vector to $\Sigma(t)$ pointing out of the fluid domain, and $\partial_t(\overline {iz_\alpha})=-i\partial_\alpha \bar z_t=-|z_\alpha|\nabla_{\bold n}\bar z_t$.}, and similar to the calculation \eqref{eq:11}, \eqref{2dat} for 2d,  an expression for $\frak a_t|N|$ in terms of $(I-\frak H_{\Sigma(t)})(\xi_{ttt}+\frak a |N| \nabla_{\bold n} \xi_t)=[\partial_t^2+a|N|\nabla_{\bold n}, \frak H_{\Sigma(t)}]\xi_t$ has been derived. 

The Taylor sign condition \eqref{taylor-s} also holds for $C^{2}$ interfaces of the water wave problem \eqref{euler} in 3d. This was proved in \cite{wu2} by an application of the Green's identity. Here we  give a heuristic argument via the maximum principle:

 Applying div to both sides of the Euler equation and using the assumption that $\text{curl}\bold v=0$ yields
$$\Delta P=-|\nabla \bold v|^2\le 0\qquad\text{in }\Omega(t).$$
Therefore from $P=0$ on the interface $\Sigma(t)$ and the maximum principle, we have $-\frac{\partial P}{\partial \bold n}\ge 0$. 

\eqref{3dquasi-linear} is then a quasilinear equation of hyperbolic type with the right hand side consisting of lower order terms. The local in time wellposedness of \eqref{3dww1} is proved in \cite{wu2} by applying energy estimates and an iterative argument to \eqref{3dquasi-linear}.

\begin{remark} An analogous derivation of the quasilinear structure can now be performed directly on
 \eqref{euler}. 
Let $D_t=\partial_t+{\bold v}\cdot\nabla$ be the material derivative, we can rewrite the first equation in  \eqref{euler} as
$$D_t\bold v+\bold k=-\nabla P.$$
Apply $D_t$ to both sides to find

\begin{equation}\label{quasi-euler}
D_t^2 \bold v+[D_t, \nabla]P=-\nabla D_t P.
\end{equation}
 From the third equation: $P=0$ on $\Sigma(t)$, we have $D_t P=0$ on $\Sigma(t)$,
so $\nabla D_t P$ points in the normal direction on the interface. We note that $$[D_t, \partial_{x_i}]P=-\partial_{x_i}{\bold v}\cdot\nabla P= -(\nabla P\cdot\nabla) {\bold v}_i=\frak a|N|\nabla_{\bold n}{\bold v}_i$$ on the fluid interface $\Sigma(t)$,\footnote{ In the second equality we used $\text{curl}\bold v=0$.} so it corresponds to the term $-\frak aN_t=\frak a|N|\nabla_{\bold n}\xi_t$ in equation \eqref{3dquasi-linear}, therefore \eqref{quasi-euler} is the counterpart of our quasilinear equation \eqref{3dquasi-linear} in the entire fluid region.
 In \cite{cl, li}, the authors used \eqref{quasi-euler} to study a more general case where the vorticity $\text{curl} \bold v$ need not be zero.\footnote{In \cite{cl, li}, the strong Taylor condition \eqref{taylor-s} is assumed to hold.}
\end{remark}

We now turn to the question of long time behavior of solutions for the water wave equation \eqref{euler} for small initial data. 

\section{Global and almost global wellposedness of the water wave equations}\label{section3}

To understand the global in time behavior of the water wave motion, we need to understand the dispersion, the nature of the nonlinearity of the water waves and their interaction.

In \cite{wu3, wu4} we studied the water wave equation \eqref{euler} in two and three space dimensions for small data, we found a nonlinear transformation for the unknowns and a nonlinear change of the coordinates, so that the transformed quantities in the new coordinate system satisfy equations containing no quadratic nonlinear terms.\footnote{When understood appropriately.} Using these canonical equations, we showed that for small, smooth and localized data of size $\epsilon$, the solution of the 2d water wave equation \eqref{euler} remain small and smooth for time $0<t<e^{c/\epsilon}$ and for similar data, the solution of the 3d water wave equation \eqref{euler} remain small and smooth for all time. 

Let's give a brief explanation of the dispersion of the water waves, the structural advantage of our canonical equations and how we found the transformations.

Let $\frak u=\bar z_t$ (or $\frak u=\xi_t$). Linearizing the quasi-linear system \eqref{2dquasilinear}-\eqref{2dat}-\eqref{al} (or \eqref{3dquasi-linear}) at  the zero solution gives
\begin{equation}\label{linearized}
\partial_t^2\frak u+|D|\frak u=F( \frak u_t, |D|\frak u),\qquad (\alpha,t)\in \mathbb R^{n-1}\times \mathbb R
\end{equation}
where $|D|=\sqrt{-\Delta}$, $\Delta$ is the Laplacian in $\mathbb R^{n-1}$ for n-dimensional water waves, $F$ consists of the nonlinear terms.  We know the dispersion relation of the  linear water wave equation 
\begin{equation}\label{linear}
\partial_t^2\frak u+|D|\frak u=0
\end{equation}
 is
\begin{equation}\label{dispersion}
\omega^2=|k|
\end{equation}
for plane wave solution $\frak u=e^{i(k\cdot\alpha+\omega t)}$, so waves of wave number $k$ travel with phase velocity
$\frac\omega{|k|}\hat k$, where $\hat k=\frac k{|k|}$, equation \eqref{linearized} is dispersive. 
For a large class $\mathfrak B$ of smooth initial data, the solution of the linear equation \eqref{linear} exists for all time and remains smooth, and its $L^\infty$ norm decays with rate $1/t^{\frac{n-1}2}$. 
  The question is  for small data in $\mathfrak B$, for how long does the solution   of the nonlinear equation  \eqref{linearized}  remain smooth.  We know nonlinear interactions can cause blow-up of the solutions at finite time. So to answer this question, we need to know for how long does the linear part of the equation \eqref{linearized} remain dominant. The weaker the nonlinear interaction, the longer the solution remains smooth. For small data, quadratic interactions are in general stronger than cubic and higher order interactions. 

To understand these assertions in qualitative terms, let's consider the following model equation with a $(p+1)$th-order nonlinearity:
\begin{equation}\label{model}
\partial_t^2\frak u+|D|\frak u=(\partial_t\frak u)^{p+1},\qquad (\alpha,t)\in \mathbb R^{n-1}\times \mathbb R.
\end{equation}
Suppose we can prove decay estimates for the solution: for $i\le s-10$, 
$$|\partial^i\partial_t\frak u(t)|_{L^\infty}\lesssim {(1+t)^{-\frac{n-1}2}} E_s(t)^{1/2},$$
where $\partial$ is some kind of derivatives, 
$$E_s(t)=\sum_{|j|\le s}\int |\partial^j\partial_t\frak u(\alpha,t)|^2+|\partial^j|D|^{1/2}\frak u(\alpha,t)|^2\,d\alpha.$$
Then we can derive energy estimates for large enough $s$:
$$\frac d{dt} E_s(t)\lesssim {(1+t)^{-\frac{(n-1)p}2}} E_s(t)^{p/2+1},$$
therefore
\begin{equation}\label{energyi}
E_s(0)^{-p/2}-E_s(T)^{-p/2}\lesssim \int_0^T{(1+t)^{-\frac{(n-1)p}2}}\,dt.
\end{equation}
Heuristically, we expect to prove existence of solutions of \eqref{model} for as long as the energy $E_s(t)$ remains finite; by \eqref{energyi}, this can be achieved if $$\int_0^T{(1+t)^{-\frac{(n-1)p}2}}\,dt\lesssim E_s(0)^{-p/2}.$$
Now if $p=1$, i.e. if the nonlinear term in \eqref{model} is quadratic, then for both $n=2$ and $n=3$, the integral $\int_0^\infty {(1+t)^{-\frac{(n-1)p}2}}\,dt=\infty$, and we would not be able to conclude solutions exist for all time for small initial data from this analysis. In fact for $p=1$ and $n=2$, the expected existence time is of order $O(\epsilon^{-2})$ for data of size $\epsilon$. If $p\ge 2$, i.e. if there is no quadratic nonlinearity, then we can expect to prove longer time existence for solutions of \eqref{model} for small initial data. In fact, we can expect for $n=2$, $p=2$  an existence time period of $[0, e^{c/\epsilon^2}]$; and for $n=3$, $p=2$ an existence time period $[0, \infty)$ for data of size $\epsilon$, when $\epsilon$ is sufficiently small.

Now for  the water wave equation \eqref{linearized}, the nonlinearity $F(\frak u_t, |D|\frak u)$ contains quadratic terms,  so appears too strong to conclude a global existence result.  
The question is whether there is another unknown $\frak v$ that satisfies an equation of the type
\begin{equation}\label{cubic1}
\partial_t^2\frak v+|D|\frak v=F_1( \frak v_t, |D|\frak v, \frak u_t, |D|\frak u)
\end{equation}
with $F_1$ containing no quadratic nonlinearities and $\|\frak v\|\approx \|\frak u\|$ in various norms $\|\cdot\|$ involved in the  analysis. The idea of finding such an unknown $\frak v$  is the so called method of normal forms, originally introduced by Poincar\'e to solve ordinary differential equations. Certainly in most cases, one should not expect such a new unknown $\frak v$ exist. For quadratic Klein-Gordon equation however, Simon and Shatah \cite{si, sh} succeeded in finding a  bilinear normal form transformation of the type
\begin{equation}\label{bilinear}
\frak v= \frak u+B( \frak u, \frak u)
\end{equation}
with $B( \frak u, \frak u)$ bilinear, canceling out the quadratic nonlinear terms  in the Klein-Gordon equation, and satisfying the norm equivalence $
\|\frak u+B(\frak u,\frak  u)\|\approx \|\frak u\|$
for small $\|\frak u\|$.  

For water wave equations \eqref{2dquasilinear}-\eqref{2dat}-\eqref{al} (or \eqref{3dquasi-linear}) however, a transformation of the type \eqref{bilinear} doesn't quite work since it has a small divisor; working with the velocity potential and the Bernoulli equation\footnote{Or equivalently the Zakharov-Craig-Sulem equation for the interface.}, a bilinear transform of the type \eqref{bilinear} has a loss of derivatives (for detailed calculations and discussions, see Appendix C of \cite{wu4}.). What we did in \cite{wu3, wu4} was to further introduce a change of the coordinates. Indeed this makes sense since when one applies a method from ODE to PDE, it is reasonable to also take into consideration the spatial variables. 

We give in the next two subsections the transforms  for the 2d and 3d water waves. The transforms are fully nonlinear. We found the transforms by first considering the 2d quasilinear equation \eqref{2dquasilinear}-\eqref{2dat}-\eqref{al}, starting with the ansatz \eqref{bilinear}, taking into considerations of the coordinate invariance of \eqref{2dquasilinear}, looking for a coordinate invariant transformation.\footnote{The bilinear transform of the type \eqref{bilinear} is not coordinate invariant.}
This entails much further efforts in understanding the bilinear transformation. What we finally arrived at is a fully nonlinear transform of the unknown function, coupled with a coordinate change. The process of finding the transforms is non-algorithmic.  
The transforms for the 3d water waves \eqref{3dquasi-linear} is obtained by naturally extending the 2d version via Clifford analysis, c.f. \cite{wu4}.

\subsection{The transformations for the 2D water waves} We give here the transformations we constructed in \cite{wu3} for the 2D water waves. Let $\Sigma(t): z=z(\alpha, t)
$, $\alpha\in\mathbb R$ be the interface  in Lagrangian coordinate $\alpha$. 

\begin{proposition}[Proposition 2.3 of \cite{wu3}]\label{prop1} Let $z=z(\alpha,t)$ be a solution of the 2d water wave system \eqref{interface-e}.
 Let 
$\Pi:=(I-\frak{H})(z-\bar z)$; let $k=k(\cdot,t):\mathbb R\to \mathbb R$ be an arbitrary diffeomorphism. Let $\zeta:=z\circ k^{-1}$, $D_t\zeta:=z_t\circ k^{-1}$ etc. be as in \eqref{2dzeta}-\eqref{2dzetaa1}. Then
\begin{equation}\label{2dcubic1}
\begin{aligned}
(&D_t^2-iA\partial_\alpha)(\Pi\circ k^{-1})
\\&=\frac{4}{\pi} \int \frac{(D_t\zeta(\alpha,t) - D_t\zeta(\beta,t))(\Im \zeta(\alpha,t) - \Im \zeta(\beta,t))}{|\zeta(\alpha,t) - \zeta(\beta,t)|^2} \partial_\beta D_t\zeta(\beta,t) d\beta \\
& \qquad + \frac{2}{\pi} \int \left(\frac{D_t\zeta(\alpha,t) - D_t\zeta(\beta,t)}{\zeta(\alpha,t) - \zeta(\beta,t)}\right)^2 \partial_\beta \Im \zeta(\beta,t) d\beta
\end{aligned}
\end{equation}
\end{proposition}
 Notice that the right hand side of equation \eqref{2dcubic1} is cubicly small if the velocity $D_t\zeta$, the height function $\Im \zeta$ of the interface and their derivatives are small, but the left hand side of \eqref{2dcubic1} still contains quadratic nonlinearities.  Naturally we ask if there is a coordinate change $k$, so that $b=k_t\circ k^{-1}$ and $A-1=(\frak ak_\alpha)\circ k^{-1}-1$ are quadratic. We need not look far, equation \eqref{b} suggests that we choose
\begin{equation}\label{2dcoordinate}
k(\alpha,t)=2\Re z (\alpha,t)-h(\alpha,t),\qquad \alpha\in\mathbb R.
\end{equation}
\begin{proposition}[Proposition 2.4 of \cite{wu3}]\label{prop2} Let $k$ be as given by \eqref{2dcoordinate}, $b=k_t\circ k^{-1}$ and $ A=(\frak a k_\alpha)\circ k^{-1}$. Let $\mathcal H$ be defined by \eqref{khilbert}. We have
\begin{equation}\label{2dab}
\begin{aligned}
 (I-\mathcal{H})b &=-[D_t\zeta,\mathcal H]\frac{\bar \zeta_\alpha-1}{\zeta_\alpha}\\
(I-\mathcal H)(A-1)&=i[D_t\zeta,\mathcal H]\frac{\partial_\alpha\overline {D_t\zeta}}{\zeta_\alpha}+i[D^2_{t}\zeta,\mathcal H]\frac{\bar \zeta_\alpha-1}{\zeta_\alpha}\end{aligned}
\end{equation}
\end{proposition}
 Propositions~\ref{prop1} and ~\ref{prop2} show that
 the quantity $\theta:=\Pi\circ k^{-1}=(I-\mathcal H)(\zeta-\bar\zeta)$ with the coordinate change $k$ given by \eqref{2dcoordinate} satisfies an equation of the type 
 \begin{equation}\label{2dcubic}
 \begin{cases}
 (\partial_t^2-i\partial_\alpha) \theta  =\mathcal G\\
 (I+\mathcal H)\theta=0
 \end{cases}
 \end{equation}
  with $\mathcal G$ containing only nonlinear terms of cubic and higher order. We make the following remarks:

\begin{remark}

 1. The transformation $I-\frak H$  and the coordinate change $k$ as given in \eqref{2dcoordinate} are  fully nonlinear in terms of the unknown function $z$ and its derivatives.\footnote{That is, $I-\frak H$ and $k$ are not finite sums of multi-linear operators of $z$ and its derivatives.} 

2. The bilinear part of the quantity $\Pi:=(I-\frak H)(z-\bar z)$ has a bounded  Fourier symbol. The coordinate change $k$ takes care of  the small divisor in the Fourier symbol of the bilinear normal form transformation. For a detailed explanation see Appendix C of \cite{wu4}.

3. 
 For $\theta$ satisfying $(I+\mathcal H)\theta=0$, $(\partial_t^2-i\partial_\alpha)\theta=
 \partial_t^2\theta+i\partial_\alpha\mathcal H\theta=(\partial_t^2+|D|)\theta+\text{quadratic}+...$. The quadratic nonlinearity  comes from $\mathcal H$, which depends nonlinearly on the unknown $\zeta$.

 4. Let $\phi$ be the velocity potential, $\psi(\alpha,t):=\phi(z(\alpha,t),t)$ be the trace of $\phi$ on the free interface. It has been shown in \cite{wu3} that the quantities $U_{k^{-1}}(I-\mathfrak H)\psi$, $U_{k^{-1}}\partial_t\Pi$ also satisfy equations of the type \eqref{2dcubic}, with their equations  given in Proposition 2.3 of \cite{wu3}.  $\partial_\alpha U_{k^{-1}}(I-\mathfrak H)\psi$, $U_{k^{-1}}\partial_t\Pi$, $\frac i2\partial_ \alpha U_{k^{-1}}\Pi$  and $\Im \partial_\alpha U_{k^{-1}}\Pi$ are near identity transforms of  the velocity $ D_t\zeta$, $2 D_t\zeta$, the acceleration $D_t^2\zeta$  and $2\Im \partial_\alpha\zeta$, see Propositions 2.5 and 2.6 of \cite{wu3}.
 
 5.  \eqref{2dcubic1}-\eqref{2dcoordinate} is used in \cite{tw} to give a rigorous justification of the NLS from the 2d water wave equation \eqref{euler}.
 
 6. The idea of changing coordinates  is  subsequently used in \cite{hi} to remove the quadratic nonlinear terms in the Burgers-Hilbert equation.
 
\end{remark}

We now give the proofs of Propositions~\ref{prop1},~\ref{prop2}.\footnote{The proofs are taken from \cite{wu3}.} 
\begin{proof} We first prove \eqref{2dcubic1}. Let $z=z(\cdot,t)$ be a solution of the water wave equation \eqref{interface-e}. Apply $(\partial_t^2-i\frak{a}\partial_\alpha)$ to $\Pi:= (I-\frak{H})(z-\bar z)$ and commute $\partial_t^2-i\frak{a}\partial_\alpha$ with $I-\frak H$ to find
\begin{equation}\label{eq:15}
 (\partial_t^2-i\frak{a}\partial_\alpha)\{(I-\frak{H})(z-\bar z)\}=(I-\frak{H})\{(\partial_t^2-i\frak{a}\partial_\alpha)(z-\bar z)\}
 -[\partial_t^2-i \frak{a}\partial_\alpha, \frak{H}](z-\bar z).
 \end{equation}
 Use \eqref{interface-e} to find $(\partial_t^2-i\frak{a}\partial_\alpha)(z-\bar z)=-2\bar z_{tt} $ then use $\bar z_t=\frak H\bar z_t$ to write $(I-\frak H)\bar z_{tt}$ as the commutator $[\partial_t, \frak H]\bar z_t$ and applying Lemma~\ref{lemma1} yields
 \begin{equation}\label{eq:16}
 (I-\frak H)\{(\partial_t^2-i\frak{a}\partial_\alpha)(z-\bar z)\}=(I-\frak H)(-2\bar z_{tt})= -2[z_t,\frak{H}]\frac{\bar z_{t\alpha}}{z_\alpha}.
 \end{equation}
 Applying Lemma~\ref{lemma1} to the second term gives
 \begin{equation}\label{eq:17}
 [\partial_t^2-i \frak{a}\partial_\alpha, \frak{H}](z-\bar z)=2[z_t,\frak{H}]\frac{z_{t\alpha}-\bar z_{t\alpha}}{z_\alpha}-\frac1{\pi\, i}\int(\frac{z_t(\alpha,t)-z_t(\beta,t)}{z(\alpha,t)-z(\beta,t)})^2\partial_\beta(z(\beta,t)-\bar z(\beta,t))\,d\beta.
\end{equation}
Subtract \eqref{eq:17} from \eqref{eq:16}. After cancelation this leaves, from \eqref{eq:15},
\begin{equation}\label{eq:10}
 (\partial_t^2-i\frak{a}\partial_\alpha)\{(I-\frak{H})(z-\bar z)\}  =-2[z_t,\frak{H}]\frac{z_{t\alpha}}{z_\alpha}+\frac1{\pi\, i}\int\big(\frac{z_t(\alpha,t)-z_t(\beta,t)}{z(\alpha,t)-z(\beta,t)}\big)^2(z-\bar z)_\beta\,d\beta.
 \end{equation}
 Because $\bar z_t$ and $\dfrac {\bar z_{t\alpha}}{z_\alpha}$ are  holomorphic, from  part 1 of Lemma~\ref{lemma2} $$[z_t,\bar {\mathfrak H}\frac1{\bar z_\alpha}]z_{t\alpha}=[z_t,\bar {\mathfrak H}]\frac{z_{t\alpha}}{\bar z_\alpha} =0.$$ In the first term of the right hand side of \eqref{eq:10} insert $[z_t,\bar {\mathfrak H}\frac1{\bar z_\alpha}]z_{t\alpha}$  to make it cubic. We have
\begin{equation}\label{309}
\begin{aligned}
 (\partial_t^2-i\frak{a}\partial_\alpha)&\{(I-\frak{H})(z-\bar z)\} 
 \\&= -2[z_t,\mathfrak H\frac1{z_\alpha}+\bar {\mathfrak H}\frac1{\bar z_\alpha}]z_{t\alpha}
 +\frac1{\pi\, i}\int\big(\frac{z_t(\alpha,t)-z_t(\beta,t)}{z(\alpha,t)-z(\beta,t)}\big)^2(z_\beta-\bar z_\beta)\,d\beta.
\end{aligned}
\end{equation}
Precomposing with $k^{-1}$ and expanding the two terms on the right gives \eqref{2dcubic1}.
\end{proof}

\begin{proof}
We prove Proposition~\ref{prop2}. We have from \eqref{2dcoordinate} $$k-\bar z=z-h.$$
Recall $h(\alpha,t)=\Phi(z(\alpha,t),t)$ where $\Phi(\cdot,t): \Omega(t)\to P_-$ is the Riemann mapping satisfying $\lim_{z\to\infty} \Phi_z(z,t)=1$. We have $$h_t=\Phi_t\circ z+(\Phi_z\circ z) z_t,\qquad h_\alpha=(\Phi_z\circ z )z_\alpha,$$
therefore
\begin{equation}\label{310}
\bar z_t-k_t=\Phi_t\circ z+(\Phi_z\circ z-1) z_t,\quad \bar z_\alpha-k_\alpha=(\Phi_z\circ z-1) z_\alpha.\end{equation}
Apply $(I-\frak H)$ to the first equality in \eqref{310}. Because  $\Phi_t$, $\Phi_z$ are holomorphic in $\Omega(t)$ with $\lim_{z\to\infty} \Phi_z(z,t)=1$, using Proposition~\ref{prop:hilbe} and   rewriting $(I-\frak{H})\{(\Phi_z\circ z-1)z_t\}$ as $[z_t,\frak{H}](\Phi_z\circ z-1)$ yields
\begin{equation}\label{312}
\begin{aligned}
-(I-\frak{H})&k_t=(I-\frak{H})(\bar z_t-k_t)\\&=(I-\frak{H})\{(\Phi_z\circ z-1)z_t\}=[z_t,\frak{H}](\Phi_z\circ z-1)=[z_t,\frak{H}]\frac{\bar z_\alpha-k_\alpha}{z_\alpha}.
\end{aligned}
\end{equation}
Precomposing with $k^{-1}$ gives the first equality of \eqref{2dab}.

 Now multiply $i\frak a$ then apply $(I-\frak H)$ to the second equality in  \eqref{310}. Using \eqref{interface-e} and the fact that $(I-\frak H)(\Phi_z\circ z-1)=0$, we  also have
$$(I-\frak{H})(i\frak a \bar z_\alpha-i\frak{a} k_\alpha)=(I-\frak{H})(i\frak a  z_\alpha (\Phi_z\circ z-1))=[z_{tt},\frak H](\Phi_z\circ z-1)=[z_{tt},\frak{H}]\frac{\bar z_\alpha-k_\alpha}{z_\alpha}.$$
Use  \eqref{interface-e} and Lemma~\ref{lemma1} to calculate 
$$ (I-\frak{H})(i\frak a \bar z_\alpha)= (I-\frak{H})(-\bar z_{tt}+i)=i-[\partial_t, \frak H]\bar z_t=i-[z_t,\frak H]\frac{\bar z_{t\alpha}}{z_\alpha}$$
so
\begin{equation}\label{313}
\begin{aligned}
-(I-\frak{H})&(i\frak{a} k_\alpha)=-(I-\frak{H})(i\frak a \bar z_\alpha)+[z_{tt},\frak H]\frac{\bar z_\alpha-k_\alpha}{z_\alpha}\\&=-i+[z_t,\frak H]\frac{\bar z_{t\alpha}}{z_\alpha}+[z_{tt},\frak{H}]\frac{\bar z_\alpha-k_\alpha}{z_\alpha}.
\end{aligned}
\end{equation}
Precomposing with $k^{-1}$ yields the second equality of \eqref{2dab}.
\end{proof}
To extend the 2d coordinate change \eqref{2dcoordinate} to 3d, we need an expression that does not rely on the Riemann mapping. Observe that for the diffeomorphism $k$ given by \eqref{2dcoordinate}, $k-\bar z= z-h$ and $z-h$ is holomorphic with $\Im(z-h)=\Im z$, so 
we can replace $z-h$ by $\frac12(I+\frak H)(I+\frak K)^{-1}(z-\bar z)$, where $\frak K=\Re\frak H$ is the double layer potential operator, and 
\begin{equation}\label{eq:12}
k=\bar z+ \frac12(I+\frak H)(I+\frak K)^{-1}(z-\bar z)
\end{equation}
modulo a real constant.\footnote{$z-h$ and $\frac12(I+\frak H)(I+\frak K)^{-1}(z-\bar z)$ are holomorphic with the same imaginary part, so the difference between them is a constant in $\mathbb R$.}
The expression \eqref{eq:12} is directly extendable to 3d.

\subsection{The transformation for the 3D water waves}

We extend the 2d transformations to 3d in the framework of  the Clifford algebra $\mathcal C(V_2)$. Besides those in subsection~\ref{3dquasi-l}, we need some additional notations. 

An element $\sigma\in \mathcal C(V_2)$ can be represented uniquely by
$\sigma=\sigma_0+\sum_{i=1}^3 \sigma_i e_i$, with $\sigma_i\in \mathbb R$ for $0\le i\le  3$.  Define  $\Re \sigma:=\sigma_0$ and call it the real part of $\sigma$. We call $\sigma$ a vector if $\sigma_0=0$. If not specified, we always assume in an expression $\sigma=\sigma_0+\sum_{i=1}^3 \sigma_i e_i$ that  $\sigma_i\in \mathbb R$ for $0\le i\le 3$. Define $\overline \sigma:=e_3\sigma e_3$, the conjugate of $\sigma$. 
We identify a point or a vector $\xi=(x_1, x_2, y)\in \mathbb R^3$ with its $\mathcal C(V_2)$ counterpart $\xi=x_1e_1+x_2e_2+ye_3$.  
 For vectors $\xi ,\ \eta \in \mathcal C(V_2)$, we know
\begin{equation}\label{vector1}
\xi\eta=-\xi\cdot\eta+\xi\times\eta,
\end{equation}
where $\xi\cdot\eta$ is the dot product, $\xi\times\eta$ the cross product. For vectors $\xi$, $\zeta$, $\eta$, $\xi(\zeta\times \eta)$ is obtained by first finding the cross product $\zeta\times \eta$, then regard it as a Clifford vector and calculating its multiplication with 
$\xi$ by the rule \eqref{product}. We write $\nabla=(\partial_{x_1},\partial_{x_2},\partial_y)$. 
 We abbreviate notations such as 
 \begin{align*} \frak H_{\Sigma}f(\alpha,\beta)&
=\iint K(\eta(\alpha',\beta') -\eta(\alpha,  \beta))\,(\eta'_{\alpha'}\times\eta'_{\beta'})\,f( \alpha',\beta'
)\,d\alpha' d\beta'
\\&:= \iint K(\eta' -\eta)\,(\eta'_{\alpha'}\times\eta'_{\beta'})\,f'\,d\alpha'd\beta':= \iint K\,N'\,f'\,d\alpha'd\beta'.
\end{align*}

As in the 2d case, $\frak H_{\Sigma}^2=I$ in $L^2$, and $\frak H_{\Sigma}1=0$.

 We  give the transformation for the 3D water wave equation \eqref{3dww1}.  Let the free interface $\Sigma(t)$ be given by $\xi= \xi(\alpha,\beta,t)=x_1(\alpha,\beta,t)e_1+x_2(\alpha,\beta,t)e_2+y(\alpha,\beta,t) e_3 $  in Lagrangian coordinates $(\alpha,\beta)$ with $N=\xi_\alpha\times\xi_\beta$ pointing out of the fluid domain $\Omega(t)$. For fixed $t$, 
let $k=k(\cdot,t)=k_1e_1+k_2e_2:\mathbb R^2\to \mathbb R^2$ be a diffeomorphism with Jacobian $J(k(t))>0$. Let $k^{-1}$ be such that $k\circ k^{-1}(\alpha, \beta, t)=\alpha e_1+\beta e_2$. Define 
\begin{equation}\label{3dzeta}
\zeta:=\xi\circ k^{-1}, \quad b:=k_t\circ k^{-1}, \quad A\circ ke_3:={\frak a}J(k)e_3:=\frak a k_\alpha\times k_\beta\end{equation}
Let $D_t:=U_k^{-1}\partial_t U_k$ be the material derivative, $ \mathcal N:=\zeta_\alpha\times\zeta_\beta$.
By  the chain rule, we know
\begin{equation}
D_t=\partial_t +b\cdot (\partial_\alpha,\partial_\beta) , \quad\quad U_k^{-1}(\frak{a}N\times \nabla)U_k=A\mathcal N\times \nabla=A(\zeta_\beta\partial_\alpha-\zeta_\alpha\partial_\beta),
\end{equation}
and $U_k^{-1}\frak H_{\Sigma(t)} U_k :=\mathcal H_{\Sigma(t)}  $, with 
\begin{equation}
\mathcal H_{\Sigma(t)}  f(\alpha,\beta,t)=\iint K(\zeta(\alpha',\beta',t)-\zeta(\alpha,\beta,t))(\zeta'_{\alpha'}\times\zeta'_{\beta'}) f(\alpha',\beta',t)\,d\alpha'\,d\beta'.
\end{equation}
We have
\begin{proposition}[Proposition 1.3 of \cite{wu4}]\label{prop:cubic1} Let $\xi=\xi(\alpha,\beta,t)$ be a solution of the 3d water wave system \eqref{3dww1}.
Let $\Pi=(I-\frak H_{\Sigma(t)})(\xi-\bar\xi)$, and for fixed $t$, $k(\cdot,t):\mathbb R^2\to \mathbb R^2$ be a diffeomorphism. We have 
\begin{equation}\label{3dcubic1}
\begin{aligned}
(&D_t^2-A\mathcal N\times \nabla)(\Pi\circ k^{-1})\\&=2\iint K(\zeta'-\zeta)\,(D_{t}\zeta-D'_t\zeta')\times
(\zeta'_{\beta'}\partial_{\alpha'}-\zeta'_{\alpha'}\partial_{\beta'})\overline{D_t'\zeta'}\,d\alpha'd\beta'\\
&-\iint K(\zeta'-\zeta)\,(D_{t}\zeta-D'_t\zeta')\times
((D'_t\zeta')_{\beta'}\partial_{\alpha'}- (D'_t\zeta')_{\alpha'}\partial_{\beta'})(\zeta'-\bar{\zeta'})\,d\alpha'd\beta'\\
&-\iint D_t K(\zeta'-\zeta)\,(D_{t}\zeta-D'_t\zeta')\times
(\zeta'_{\beta'}\partial_{\alpha'}-\zeta'_{\alpha'}\partial_{\beta'})(\zeta'-\bar{\zeta'})\,d\alpha'd\beta'
\end{aligned}
\end{equation}
\end{proposition}
Observe that the second and third terms in the right hand side of \eqref{3dcubic1} are cubicly small provided the velocity $D_t\zeta$ and the steepness of the height function $\partial_\alpha(\zeta-\bar\zeta)$, $\partial_\beta(\zeta-\bar\zeta)$ are small, while the first term appears to be only quadratically small. Unlike the 2D case, multiplications of Clifford analytic functions are not necessarily analytic, so we cannot reduce the first term at the right hand side of equation \eqref{3dcubic1}  into a cubic form. However we note that 
the first term is  almost  analytic in the fluid domain $\Omega(t)$, while 
   the left hand side of \eqref{3dcubic1} is almost analytic in the air region. 
    The orthogonality of the projections 
$\frac12(I-\mathcal H_{\Sigma(t)} )$ and $\frac12(I+\mathcal H_{\Sigma(t)} )$ allows us to reduce the first term into a cubic in energy estimates, see \cite{wu4}.

Now the left hand side of \eqref{3dcubic1} still contains quadratic terms. As in the 2D case,  we resolve this difficulty  by choosing an appropriate coordinate change $k$. Let
\begin{equation}\label{k}
k=k(\alpha,\beta,t)=\xi(\alpha,\beta,t)-(I+\frak H_{\Sigma(t)} )y(\alpha,\beta,t)e_3+\frak K_{\Sigma(t)}  y (\alpha,\beta,t) e_3
\end{equation}
Here $\frak K_{\Sigma(t)} =\Re \frak H_{\Sigma(t)} $: 
\begin{equation}\label{doublelayer1}
\frak K_{\Sigma(t)}  f(\alpha,\beta, t)=-\iint K(\xi(\alpha',\beta',t)-\xi(\alpha,\beta, t))\cdot N' f(\alpha', \beta',t)\,d\alpha'\,d\beta'
\end{equation}
is the double layer potential operator. It is clear that the $e_3$ component of $k$ as defined in \eqref{k} is zero. In addition, the real part of $k$ is also zero.
This is because
\begin{equation*}
\begin{aligned}
\iint K(\xi'-\xi) \times &(\xi'_{\alpha'}\times \xi'_{\beta'}) y' e_3\,d\alpha'\,d\beta'\\&=\iint (\xi'_{\alpha'} \xi'_{\beta'}\cdot K-\xi'_{\beta'}\xi'_{\alpha'}\cdot K)y'e_3\,d\alpha'\,d\beta' \\
&= -2\iint (\xi'_{\alpha'} \partial_{\beta'}\Gamma(\xi'-\xi)-\xi'_{\beta'}\partial_{\alpha'}\Gamma(\xi'-\xi))y'e_3\,d\alpha'\,d\beta'\\&
=2 \iint \Gamma(\xi'-\xi)(\xi'_{\alpha'} y_{\beta'}-\xi'_{\beta'}y_{\alpha'})e_3\,d\alpha'\,d\beta'
\\&=2\iint \Gamma(\xi'-\xi)(N_1'e_1+N_2'e_2)\,d\alpha'\,d\beta'
\end{aligned}
\end{equation*}
So
\begin{equation}\label{hilbertz}
\frak H_{\Sigma(t)}  y e_3=\frak K_{\Sigma(t)}  y e_3+2\iint \Gamma(\xi'-\xi)(N_1'e_1+N_2'e_2)\,d\alpha'\,d\beta'
\end{equation}
This shows that the mapping $k$ defined in \eqref{k} has only the $e_1$ and $e_2$ components $k=(k_1, k_2)=k_1 e_1+k_2 e_2$. If $\Sigma(t)$ is a graph with small steepness, i.e. if $y_\alpha$ and $y_\beta$ are small, then the Jacobian of $k=k(\cdot,t)$: $J(k)=J(k(t))=\partial_\alpha k_1\partial_\beta k_2-
\partial_\alpha k_2\partial_\beta k_1>0$ and $k(\cdot,t) :\mathbb R^2\to \mathbb R^2$ defines a valid coordinate change (c.f. \cite{wu4}). 

The following proposition shows that if $k$ is as given in \eqref{k}, then  $b$ and $A-1$ are  quadratic. Let 
\begin{equation}\label{zeta1}
\mathcal K_{\Sigma(t)}: =\Re \mathcal H_{\Sigma(t)} =:U_k^{-1}\frak K_{\Sigma(t)} U_k, \quad P:=\alpha e_1+\beta e_2,\quad\text{and}\quad \zeta:=P+\lambda.
\end{equation}
\begin{proposition}[Proposition 1.4 of \cite{wu4}]\label{propAb} Let $k$ be as given in \eqref{k}. Let $b=k_t\circ k^{-1}$ and  $A\circ k={\frak a}J(k)$. We have 
\begin{equation}\label{3dab}
\begin{aligned}
b&= \frac12(\mathcal H_{\Sigma(t)} -\overline{\mathcal H_{\Sigma(t)} })\overline {D_t\zeta}-\frac12[D_t,\mathcal H_{\Sigma(t)} -\mathcal K_{\Sigma(t)} ](\zeta-\bar\zeta)+ \frac12\mathcal K_{\Sigma(t)}  (D_t\zeta-\overline{D_t\zeta})\\
(A-1) e_3 &=\frac12(-\mathcal H_{\Sigma(t)} +\overline{\mathcal H_{\Sigma(t)} })\overline {D_t^2\zeta}+
\frac12 ([D_t, \mathcal H_{\Sigma(t)} ]D_t\zeta-\overline{[D_t, \mathcal H_{\Sigma(t)} ]D_t\zeta})\\&
+\frac12[A \mathcal N\times \nabla, \mathcal H_{\Sigma(t)} ](\zeta-\bar\zeta)
-\frac12A \zeta_\beta\times (\partial_\alpha\mathcal K_{\Sigma(t)} (\zeta-\bar\zeta))\\&+\frac12 A\zeta_\alpha\times (\partial_\beta\mathcal K_{\Sigma(t)}  (\zeta-\bar\zeta)  )+A\partial_\alpha\lambda\times \partial_\beta\lambda
\end{aligned}
\end{equation}
Here $\overline{\mathcal H_{\Sigma(t)} }=e_3\mathcal H_{\Sigma(t)}  e_3$.
\end{proposition}

Let $\chi=\Pi\circ k^{-1}$ with $k$ be given by \eqref{k}. The left hand side of equation \eqref{3dcubic1} is
\begin{equation*}(\partial_{t}^2-e_2\partial_\alpha+e_1\partial_\beta)\chi-\partial_\beta\lambda\partial_\alpha\chi+\partial_\alpha\lambda\partial_\beta\chi
+\text{cubic and higher order terms}
\end{equation*}
The quadratic term $\partial_\beta\lambda\partial_\alpha\chi-\partial_\alpha\lambda\partial_\beta\chi$
is new in 3d. Observe that this is one of the null forms studied in \cite{kl4}.  It is also null for our equation and can be written as the factor $1/t$ times a quadratic expression involving some "invariant vector fields" for $\partial_{t}^2-e_2\partial_\alpha+e_1\partial_\beta$, see \cite{wu4}. Therefore this term is cubic in nature and  equation \eqref{3dcubic1} is of the type "linear + cubic and higher order perturbations".

We refer the reader to \cite{wu4} for the proof of Propositions~\ref{prop:cubic1} and ~\ref{propAb}.

We remark that the 3d transformations is recently used in \cite{to} to give a rigorous justification of the modulation approximation for the 3d water wave equation \eqref{euler}. Besides the 3d transforms, \cite{to} uses the method of normal form to handle the quadratic term $\partial_\beta\lambda\partial_\alpha\chi-\partial_\alpha\lambda\partial_\beta\chi$  since this term is truly quadratic in the modulation regime.

\subsection{Global in time behavior of solutions for the 2d and 3d water waves} In \cite{wu3, wu4}  equations \eqref{2dcubic1} and \eqref{3dcubic1} together with the coordinate changes \eqref{2dcoordinate} and \eqref{k} are used to prove the almost global wellposedness of the 2d and global wellposedness of the 3d water wave equations \eqref{euler} for small, smooth and localized initial data. The basic idea is what we illustrated with the model equation \eqref{model}; it is made rigorous by the method of invariant vector fields.  This involves constructing invariant vector fields for the operator $\partial_t^2-e_2\partial_\alpha+e_1\partial_\beta$ (the invariant vector fields for $\partial_t^2-i\partial_\alpha$ for the 2d case is available due to the well studied Schr\"odinger operator $i\partial_t-\partial_x^2$), proving generalized Sobolev inequalities that give  $L^2\to L^\infty$ estimates with the decay rate  $1/{t^{1/2}}$ for the 2d and $1/t$ for the 3d water waves, using equations \eqref{2dcubic1}-\eqref{2dcoordinate} and \eqref{3dcubic1}-\eqref{k}  to show that properly constructed energies that involve invariant vector fields 
remain bounded for the time period $[0, e^{c/\epsilon}]$ for the 2d  and for all time for the 3d water waves for data of size $O(\epsilon)$. 
 The projection $\frac12(I-\frak H)$ is used  in various ways to project away  "quadratic noises" in the course of deriving the energy estimates.  We remark that it is more natural to treat $D_t^2-iA\partial_\alpha$ and  $D_t^2-A\mathcal N\times \nabla $ as the main operators for the 2d and 3d water wave equations than treating them as the perturbations of the linear operators $\partial_t^2-i\partial_\alpha$ and $\partial_t^2-e_2\partial_\alpha+e_1\partial_\beta$. 
The almost global well-posedness for the 2d and global well-posedness for the 3d water wave equations follow from the local well-posedness results, the uniform boundedness of the energies  and continuity arguments. For details of the proofs see \cite{wu3, wu4}. We state the results.

Let $|D|=\sqrt{-\Delta}$, $H^{s}(\mathbb R^{n-1})=\{ f\, |\, (I+|D|)^{s}f\in L^2(\mathbb R^{n-1})\}$, with $\|f\|_{H^{s}(\mathbb R^{n-1})}=\|(I+|D|)^{s}f\|_{L^2(\mathbb R^{n-1})}.$

\subsection*{2d water waves}
Let $s\ge 12$, $\max\{[\frac s2]+3, 11\}\le l \le s-1$. Assume  
\begin{equation}\label{2dic}
\begin{aligned}
z(\alpha,0)=(\alpha, y(\alpha)), \quad  z_t(\alpha,0)&=\frak u(\alpha),\quad z_{tt}(\alpha,0)=\frak w(\alpha) \quad \alpha\in \mathbb R, \\
 \bold v(z,0)&=\bold g(z),\quad z\in \Omega(0)
\end{aligned}
\end{equation}
and the data in \eqref{2dic} satisfies the 2d water wave system \eqref{interface-e}. In particular $\bar{\bold g}$ is a holomorphic function in the initial fluid domain $\Omega(0)$ and $\bold g(z(\alpha,0))=\frak u(\alpha)$. Let $\Gamma=\partial_\alpha, \alpha\partial_\alpha$. Assume that
$$\sum_{|j|\le s-1}(\|\Gamma^j y_\alpha\|_{H^{1/2}(\mathbb R)}+\|\Gamma^j \frak u\|_{H^{3/2}(\mathbb R)}+\|\Gamma^j\frak w\|_{{H^1}(\mathbb R)})<\infty.$$
Let 
$$\epsilon=\sum_{|j|\le l} (\|\Gamma^j y\|_{H^1(\mathbb R)}+\|\Gamma^j\frak u\|_{H^1(\mathbb R)})+\sum_{j\le l-2}\|(z\partial_z)^j\bar {\bold g}\|_{L^2(\Omega(0))}.$$
\begin{theorem}[2d Theorem, c.f.\cite{wu3}]\label{2dtheorem} There exist $\epsilon_0$ and $c>0$, such that for $\epsilon<\epsilon_0$, the initial value problem \eqref{interface-e}-\eqref{2dic} has a unique classical solution for the time period $[0, e^{c/\epsilon}]$. During this time, the interface is a graph, the solution is as regular as the initial data and remains small. Moreover the $L^\infty$ norm of the steepness  $\partial_\alpha(z-\bar z)$, the velocity $z_t$ and acceleration $z_{tt}$ decay at rate $1/{t^{1/2}}$.
\end{theorem}

\subsection*{3d water waves} 
Let $s\ge 27$, $\max\{[\frac s2]+1, 17\}\le l\le s-10$. Assume that initially  
\begin{equation}\label{ic}
\xi(\alpha,\beta,0)=(\alpha,\beta, y^0(\alpha,\beta)),\quad \xi_t(\alpha,\beta,0)=\frak u^0(\alpha,\beta),\quad \xi_{tt}(\alpha,\beta,0)=\frak w^0(\alpha,\beta),
\end{equation}
 and the  data in \eqref{ic} satisfies the 3d water wave system \eqref{3dww1}.  Let $\Gamma=\partial_\alpha,\,\partial_\beta,\, \alpha\partial_\alpha+\beta\partial_\beta,\, \alpha\partial_\beta-\beta\partial_\alpha$. Assume that
\begin{equation}
\sum_{|j|\le s-1\atop\partial=\partial_\alpha,\partial_\beta}\|\Gamma^j |D|^{1/2} y^0\|_{L^2(\mathbb R^2)}+\|\Gamma^j\partial y^0\|_{H^{1/2}(\mathbb R^2)}+\|\Gamma^j \frak u^0\|_{H^{3/2}(\mathbb R^2)}+\|\Gamma^j\frak w^0\|_{H^1(\mathbb R^2)}<\infty
\end{equation}
Let
\begin{equation}\label{assumemain}
\epsilon=\sum_{|j|\le l+3\atop\partial=\partial_\alpha,\partial_\beta}\|\Gamma^j |D|^{1/2} y^0\|_{L^2(\mathbb R^2)}+\|\Gamma^j\partial y^0\|_{L^2(\mathbb R^2)}+\|\Gamma^j \frak u^0\|_{H^{1/2}(\mathbb R^2)}+\|\Gamma^j\frak w^0\|_{L^2(\mathbb R^2)}.
\end{equation}
\begin{theorem}[3d Theorem, c.f. \cite{wu4}]\label{3dtheorem}
There exists $\epsilon_0>0$, such that for $0< \epsilon\le \epsilon_0$, the initial value problem \eqref{3dww1}--\eqref{ic} has a unique classical solution globally in time. For each time $0\le t<\infty$,  the interface is a graph, the solution has the same regularity as the initial data and remains small. Moreover the $L^\infty$ norm of the steepness and the acceleration on the interface, the derivative of the velocity on the interface decay at  rate $\frac 1t$. 

\end{theorem}

\begin{remark} 1. The existence time in Theorem~\ref{2dtheorem} is extended to global in \cite{ip, ad, it} by further understanding the nature of the cubic nonlinearities in equation \eqref{euler}. 

2. In \cite{wu3} a quick dispersive estimate was proved by the vector field method. An in-depth analysis on the dispersion of the linear water wave operator $\partial_t^2+|D|$ is performed in \cite{be}; a threshold is found so that when the amount of small frequency waves in the initial data is below the threshold, the solution of the linear water wave equation \eqref{linear} decays with rate $1/{\sqrt{t}}$, while above the threshold, there is a growth factor in the linear solution. 
Consequences in the nonlinear setting remain to be understood.
\end{remark}

{\bf{Acknowledgement}}: The author thanks Jeffrey Rauch and Shuang Miao for carefully reading through the manuscript and for their remarks and suggestions.

\begin{appendix}
\section{Basic analysis preparations}
In this section we present some inequalities and identities on $\mathbb R$ that are used to guide the derivation of the quasilinear structure of the 2d water waves. 
Corresponding inequalities and identities are available in all dimensions $\mathbb R^d$, we refer the reader to \cite{wu2, wu4} for those.

Let $H\in C^1(\mathbb R; \mathbb R^l), \ A_i\in C^1(\mathbb R)$, $i=1,\dots m$, $F\in C^\infty(\mathbb R^l)$. Define
\begin{equation}\label{3.15}
C_1(H, A, f)(x)=p.v.\int F(\frac{H(x)-H(y)}{x-y})\frac{\Pi_{i=1}^m(A_i(x)-A_i(y))}{(x-y)^{m+1}}f(y)\,dy.
\end{equation}

\begin{proposition}\label{A1} There exist constants $c_1=c_1(F, \|H'\|_{L^\infty}) $, $c_2=c_2(F, \|H'\|_{L^\infty})$, such that 

1. For any $f\in L^2,\ A_i'\in L^\infty, \ 1\le i\le m, $
\begin{equation}\label{3.16}
\|C_1(H, A, f)\|_{L^2}\le c_1\|A_1'\|_{L^\infty}\dots\|A_m'\|_{L^\infty}\|f\|_{L^2}. 
\end{equation}
2. For any $ f\in L^\infty, \ A_i'\in L^\infty, \ 2\le i\le m,\ A_1'\in L^2$, 
\begin{equation}\label{3.17}
\|C_1(H, A, f)\|_{L^2}\le c_2\|A_1'\|_{L^2}\|A'_2\|_{L^\infty}\dots\|A_m'\|_{L^\infty}\|f\|_{L^\infty}.
\end{equation}
\end{proposition}
\eqref{3.16} is a result of Coifman, McIntosh and Meyer \cite{cmm}. \eqref{3.17} is a consequence of the Tb Theorem, a proof  is given in \cite{wu3}.

Let $H$, $A_i$, $F$ satisfy the same assumptions as in \eqref{3.15}. Define
\begin{equation}\label{3.19}
C_2(H, A, f)(x)=\int F(\frac{H(x)-H(y)}{x-y})\frac{\Pi_{i=1}^m(A_i(x)-A_i(y))}{(x-y)^m}\partial_y f(y)\,dy.
\end{equation}
We have the following inequalities.
\begin{proposition}\label{A2} There exist constants $c_3=c_3(F, \|H'\|_{L^\infty}) $, $c_4=c_4(F, \|H'\|_{L^\infty})$, such that 

1. For any $f\in L^2,\ A_i'\in L^\infty, \ 1\le i\le m, $
\begin{equation}\label{3.20}
\|C_2(H, A, f)\|_{L^2}\le c_3\|A_1'\|_{L^\infty}\dots\|A_m'\|_{L^\infty}\|f\|_{L^2}.
\end{equation}

2. For any $ f\in L^\infty, \ A_i'\in L^\infty, \ 2\le i\le m,\ A_1'\in L^2$,
\begin{equation}\label{3.21}
\|C_2(H, A, f)\|_{L^2}\le c_4\|A_1'\|_{L^2}\|A'_2\|_{L^\infty}\dots\|A_m'\|_{L^\infty}\|f\|_{L^\infty}.\end{equation}
\end{proposition}

Using integration by parts, the operator $C_2(H, A, f)$ can be easily converted into a sum of operators of the form $C_1(H,A,f)$. \eqref{3.20} and \eqref{3.21} follow from \eqref{3.16} and \eqref{3.17}.  

 The following identities are useful to compute the derivatives of the integral operators.

Let $${\bold K}f(x,t)=p.v.\int K(x,y;t)f(y,t)\,dy$$
where  either
$K$ or $(x-y)K(x,y;t)$ is continuous and bounded, and $K$ is  smooth away from the diagonal $\Delta=\{(x,y)\,|\,x=y\}$. 
We have for $f\in C^1(R^{1+1})$ vanishing  as $|x|\to \infty$, 
\begin{equation}\label{kd}
\begin{aligned} 
&[\partial_t, {\bold K}]f(x,t)=\int\partial_t K(x,y;t)f(y,t)\,dy
\\&[\partial_x, {\bold K}]f(x,t)=\int(\partial_x+\partial_y)K(x,y;t)f(y,t)\,dy
\end{aligned}
\end{equation}

The first identity in \eqref{kd} is straightforward,  the second is obtained by integration by parts. 

Using \eqref{kd}, $\partial_{\alpha'}^s[f,\mathbb H]\partial_{\alpha'}g$ equals to the sum of $[\partial_{\alpha'}^s f,\mathbb H]\partial_{\alpha'}g$,\footnote{Expand $[\partial_{\alpha'}^s f,\mathbb H]\partial_{\alpha'}g= \partial_{\alpha'}^s f\mathbb H(\partial_{\alpha'}g)-\mathbb H(\partial_{\alpha'}^s f\partial_{\alpha'}g)  $ and estimate term by term.} $[f,\mathbb H]\partial_{\alpha'}^{s+1}g$ and some intermediate terms. By Propositions~\ref{A1} and~\ref{A2} and the Sobolev embedding, $[f,\mathbb H]\partial_{\alpha'}g$ has the same regularity as $f$ and $g$. 

When the Riemann mapping is not used,  typically we work in the regime where the interface is chord-arc, that is
there are constants $\mu_1>0$ and $\mu_2>0$, such that
$$\mu_1|\alpha-\beta|\le |z(\alpha,t)-z(\beta,t)|\le \mu_2|\alpha-\beta|,\qquad \forall  \alpha,\beta\in\mathbb R.$$
By \eqref{kd}, derivatives of $[f,\frak H]\frac{\partial_\alpha g}{z_\alpha}$ are of types \eqref{3.15} and \eqref{3.19}, where in \eqref{3.15} or \eqref{3.19}  $F\in C^\infty_0(\mathbb C)$ is chosen to be $F(z)=\frac1{z^m}$ for $\mu_1\le |z|\le \mu_2$ and $H(\alpha)=z(\alpha,t)$.  An application of Proposition~\ref{A1} and~\ref{A2} and Sobolev embedding  shows that $[f,\frak H]\frac{\partial_\alpha g}{z_\alpha}$  has the same regularity as $f$, $g$ and $z$.

\section{Commutator identities}
We need the following identities in our derivation of the quasilinear structure \eqref{eq:11} and 
Propositions~\ref{prop1} and~\ref{prop2}.  Let $z=z(\cdot,t)$ define a non-selfintersecting curve, and $\frak H$ be the Hilbert transform as defined in \eqref{hilbert-t}.
\begin{lemma}[Lemma 2.1 of \cite{wu3}]\label{lemma1} Assume that $z_t, z_\alpha-1\in C^1([0, T], H^1(\mathbb R))$, $f\in C^1(\mathbb R\times (0, T))$ satisfies $f_\alpha(\alpha,t)\to 0$, as $|\alpha|\to \infty$. We have
\begin{equation}\label{3.11}
\begin{aligned} \ 
[\partial_t, \frak{H}]f&=[z_t,\frak{H}]\frac{f_\alpha}{z_\alpha}\\
[\partial_t^2, \frak{H}]f&=[z_{tt},\frak{H}]\frac{f_\alpha}{z_\alpha}+2[z_t,\frak{H}]\frac{f_{t\alpha}}{z_\alpha}-\frac1{\pi\, i}\int\big(\frac{z_t(\alpha,t)-z_t(\beta,t)}{z(\alpha,t)-z(\beta,t)}\big)^2f_\beta\,d\beta\\
[\frak{a}\partial_\alpha,\frak{H}]f&=[\frak{a}z_\alpha,\frak{H}]\frac{f_\alpha}{z_\alpha}, \quad  \partial_\alpha\frak{H}f=z_\alpha\frak{H}\frac{f_\alpha}{z_\alpha}\\
[\partial_t^2-i \frak{a}\partial_\alpha, \frak{H}]f&=2[z_t,\frak{H}]\frac{f_{t\alpha}}{z_\alpha}-\frac1{\pi\, i}\int\big(\frac{z_t(\alpha,t)-z_t(\beta,t)}{z(\alpha,t)-z(\beta,t)}\big)^2f_\beta\,d\beta\\
(I-\frak{H})(-i{\frak a}_t \bar{z}_\alpha)&=2[z_{tt}, \frak {H}]\frac{\bar z_{t\alpha}}{z_\alpha}+2[z_t, \frak {H}]\frac{\bar z_{tt\alpha}}{z_\alpha} -\frac 1{\pi  i}\int(\frac{z_t(\alpha, t)-z_t(\beta, t)}{z(\alpha,t)-z(\beta,t)})^2\bar z_{t\beta}\, d\beta
\end{aligned}
\end{equation}
\end{lemma}

The proof of Lemma~\ref{lemma1} is straightforward by integration by parts. We omit the proof.   

Let $\Omega\subset \mathbb C$ be a domain with boundary $\Sigma: z=z(\alpha), \alpha\in I$ oriented clockwise. Let $\frak H$ be defined by \eqref{hilbert-t}. 

\begin{lemma}[Lemma 2.2 of \cite{wu3}]\label{lemma2} 1. If $f=\frak H f$, $g=\frak Hg$, then $[f, \frak H]g=0$.

2. For any $f, g\in L^2$, we have $[f, \frak H]\frak H g=-[\frak H f,\frak H]g$.

\end{lemma}

The first statement follows from the fact that the product of holomorphic functions is holomorphic. Observe that the first statement also holds for $f=-\frak Hf$ and $g=-\frak Hg$, the second statement follows from the first by applying the first identity to $(I\pm \frak H)f$ and $(I\pm\frak H)g$.

We do not give the commutator identities for 3d, but refer the reader to \cite{wu2, wu4} for details. 

\end{appendix}

\end{document}